\documentclass{amsart}

\usepackage{amsthm,amsfonts,amsmath,amssymb,latexsym,epsfig}
\usepackage{upref,amssymb,eucal,ae,enumitem}
\usepackage[all,cmtip]{xy}

\newtheorem{theorem}{Theorem}
\newtheorem{lemma}[theorem]{Lemma}

\newtheorem{corollary}[theorem]{Corollary}
\newenvironment{example}{\medskip \refstepcounter{theorem}
\noindent  {\bf Example \thetheorem}.\rm}{\,}
\newenvironment{remark}{\medskip \refstepcounter{theorem}
\noindent  {\bf Remark \thetheorem}.\rm}{\,}

\def\hok{\mbox{}\begin{picture}(10,10)\put(1,0){\line(1,0){7}}
\put(8,0){\line(0,1){7}}\end{picture}\mbox{}}

\def\<{\langle}
\def\>{\rangle}

\def\tn{\tilde{n}}
\def\tg{\tilde{g}}
\def\mb#1{{\mathbb #1}}
\def\mc#1{{\mathcal #1}}
\def\mf#1{{\mathfrak #1}}
\def\BOne{{\mathchoice {\rm 1\mskip-4mu l} {\rm 1\mskip-4mu l}
                          {\rm 1\mskip-4.5mu l} {\rm 1\mskip-5mu l}}}

\begin{document}

\title[Vector field cycles in the tangent bundle]
{Vector field cycles in the tangent bundle}
\author{Santiago R. Simanca}
\email{srsimanca@gmail.com}

\begin{abstract}
Given a closed Riemannian manifold $(M^m,g)$ and a vector field $v$ on $M$, we
form the Sasaki metric $g_S$ on $TM$, and restrict it to the image of the cross
section map of $M$ into $TM$ defined by $v$, whose pull back to $M$ defines a
new metric $g(v)$ on $M$. We then view the cross section as an isometric 
embedding $f_{g(v)}: (M,g(v))\rightarrow (TM,g_S)$, which when $\| v\|_g=1$, 
ranges into the unit sphere bundle $(S^1(TM),g_S)$. $v$ is 
minimal or minimal unit if these embeddings have null mean curvature vectors, 
conditions that occur if, and only if, $v$ is in the kernel or $v$ is an 
eigenvector, respectively, of a first order perturbation of a weighted rough 
Laplacian, the 
weights and perturbation determined by the covariant derivatives 
$\nabla^g_{e}v$ along unit directions $e$ in suitable orthonormal frames that 
include $v$ when $\| v\|_g=1$, and the curvature tensor of $g$. 
A minimal unit field must be Killing, and other than parallel 
fields, $v=0$ is the only minimal one. We characterize the minimal 
unit vector fields on the standard sphere $(\mb{S}^{2n+1},g)
\hookrightarrow (\mb{R}^{2n+2},\| \,  \|^2)$ as those defining  
contact strictly pseudoconvex CR structures whose Levi form and sign 
are determined by $g$ and the orientation.  
If $\Theta_{f_{g(v)}}(M)$ and $\Phi_{f_{g(v)}}(M)$ 
are the total exterior scalar
curvature and squared $L^2$ norm of the mean curvature vector functionals, 
and $m>2$, a canonical cycle $f_{g(v)}(M)$ (in the 
homology class $[f_{g(v)}(M)]$) is a critical point of the 
functional $(m/m-1) \Theta_{f_{g(v)}}(M) +\Phi_{f_{g(v)}}(M)$ under (volume 
preserving) conformal deformations, notions conveniently defined 
also when $m\leq 2$.
The zero section of $TM$ is a canonical cycle if, and only if, 
the scalar curvature 
of $g$ is constant. We describe some examples of these vector fields
and cycles, and analyze their deformations under dilations of the field.
\end{abstract}


\maketitle

\section{Introduction}
Suppose that $v \in C^{\infty}(M;TM)$ is a vector field on a closed Riemannian
manifold $(M,g)$. Then we have the natural embedding
$$
\begin{array}{ccc}
M & \stackrel{f(v)}{\rightarrow} & TM \\
p & \rightarrow & (p,v(p))  
\end{array} \, , 
$$ 
and using the Sasaki metric $g_S$ on $TM$, 
if $\iota_v: f(v)(M)\hookrightarrow TM$, then $g(v)^{TM}:=\iota_v^*g_S$ is a 
metric on $f(v)(M)$ whose  pull back $f(v)^*g(v)^{TM}$ to $M$ 
defines a new metric $g(v)$ on $M$. We set  $\mu^{TM}_{g(v)}(v):=
\mu_{{g(v)}}(M)=\mu_{{g(v)^{TM}}}(f(v)(M))$, and call it the volume of $v$.   
We obtain a family of isometric embedding mappings
\begin{equation} \label{egv}
C^{\infty}(M;TM)\ni v \mapsto
(M,g(v)) \stackrel{f_{g(v)}}{\rightarrow} (TM,g_S)\, .
\end{equation}
We begin by computing the mean curvature vector $H^{TM,g_S}_{f_{g(v)}}$ of an
$f_{g(v)}$, and show that it is null 
if for any orthonormal frame of eigenvectors of $(\nabla^{g}v)^t \nabla^g v$
(where $\nabla^gv$ is viewed as a $(1,1)$ tensor),   
$v$ is in the kernel of a first order 
perturbation of a weighted rough Laplacian, with the weights determined in a 
certain way by 
the covariant derivatives $\nabla^g_e v$ along unit directions in
the frame, and perturbation by these fields and curvature tensor $R^g$ of $g$, 
respectively. It follows that $f_{g(v)}$ is a minimal embedding if, and only 
if, the vector field $v$ is parallel or null, in which case,
$\mu^{TM}_{g(v)}(v)=\mu_g(M)$ is the absolute minimum of the volume 
functional of $v$. 

Let us suppose now that $(M,g)$ carries nonvanishing vector fields. If 
$v \in C^{\infty}(M;S^1(TM))$ is a unit vector field, then $f_{g(v)}$ embeds
$M$ into the unit sphere bundle $S^1(TM) \hookrightarrow TM$, and in this
context, we obtain a family of isometric embedding mappings
\begin{equation} \label{egvU}
C^{\infty}(M;S^1(TM))\ni v \mapsto
(M,g(v)) \stackrel{f_{g(v)}}{\rightarrow} (S^1(TM),g_S)\, ,
\end{equation}
where the metric on the background sphere bundle is that induced by the 
metric $g_S$ on $TM$, which we denote by $g_S$ also. We let
$\mu^{S^1(TM)}_{g(v)}(v)=\mu_{g(v)}(M)$ be the volume functional of the unit
$v$. 

In the 1980s, thinking about standard odd dimensional spheres, 
Gluck and Ziller \cite{glzi} asked if the geometrically best unit vector 
fields on them would always realize the minimum volume. 
By using the calibrations of Federer and Harvey-Lawson, they 
proved that the only unit vector fields of minimal volume on $\mb{S}^3$,
up to conjugation, are Hopf vector fields, and by exploiting the 
method for doing so, concluded that these fields are the unique absolute 
volume minimizers in their (common) homology class. Soon after, 
Johnson \cite{john} proved that in dimension $5$ or above, the Hopf vector 
fields are unstable critical points of the volume functional, while Pedersen 
\cite{pede} proved that there exist unit vector fields that converge to a 
vector field with one singularity while the volume gets arbitrarily close to 
the minimum. 

Ever since, the various questions raised by Gluck and Ziller have generated a 
tremendous amount of research, evidence of which can be found in the 
exhaustive bibliography of the recent monograph of Gil-Medrano \cite{ol}, 
who proposed the study of the critical points of the volume  
functional $\mu^{S^1(TM)}_{g(v)}(v)$ as a way of making inroads into addresing 
them. In particular, the critical point equation under variations of $v$, 
and its second variation at a critical point, were derived 
in \cite{olel} and \cite{bogi1,olel2}, respectively, both as suitable 
equations in $(TM,g_S)$ ported to $(M,g)$ using the embedding; these 
were then used to characterize the scaled Hopf vector fields on odd dimensional
spheres as minimal (in the scaled sense), and to determine their stability as 
a function of the scale. About this time,  the lower bound of the volume 
functional of unit fields on these spheres of Gluck and Ziller was generalized 
to  spaces of constant curvature \cite{brcn}, and shown that in dimensions
greater than three, the bound cannot be achieved, though it is so by a field 
with two singularities. But in the existing literature, there seems 
to be a dearth of examples other than minimal unit fields that are Killing,
and even the problem of characterizing the minimal unit vector fields on the 
standard sphere has remained open. 

Guided by the approach used in the earlier problem, which directly ties the 
volume functional and the isometric embedding, we compute the mean 
curvature vector $H^{S^1(TM),g_S}_{f_{g(v)}}$ of an $f_{g(v)}$ in (\ref{egvU}),
and show that it is null if for any orthonormal frame consisting of
$v$ and eigenvectors of $(\nabla^{g}v)^t \nabla^g v$, $v$ is now an eigenvector 
of the same first order perturbation of the weighted rough Laplacian of 
earlier, of weights and perturbation determined by the local fields 
$\nabla^g_e v$ along unit directions $e$ in the frame, and curvature tensor 
$R^g$, with the square of the eigenvalue now equaling the squared pointwise 
norm of the mean curvature vector of $f_{g(v)}$ viewed as an embedding into 
$(TM,g_S)$. It follows that 
a $v$ that is both minimal and minimal unit must be parallel, and by exploiting
the freedom of choice in the frame in which we express the minimality condition,
we prove that a minimal unit vector field must be Killing, 
filling one of the voids in the literature alluded to earlier. This
allows us to exploit the properties of the foliations these fields define, 
perhaps in the spirit of Sullivan's \cite{sull,de}, where the nonintegrable
hyperplane of transversals play a protagonist role. These properties 
are the key into our characterization of the minimal unit vector fields on
a standard odd dimensional sphere, which we prove are in 1-to-1 correspondence
with contact strictly pseudoconvex CR structures on the sphere whose
Levi form and sign are determined by $g$ and the orientation.  

We organize the paper as follows: In \S2, which serves to set up notation, 
we discuss the specifics of connections on a vector bundle that we need,
as well as those of the Sasaki metric on $TM$. In  \S3, we derive the 
equations for the minimality of isometric embeddings of
type (\ref{egv}) and (\ref{egvU}), respectively, and prove that the vector
field that defines a critical point of the latter type must be a Killing unit
vector field. By exploiting the local nature of the critical point equation, we
extend here the notion of minimal unit to vector fields with singulaties 
satisfying a mild condition. We discuss also the extrinsic quantities of any 
embedded cycle $f_{g(v)}(M)$, and combine them appropriately to define 
suitable functionals whose critical points single out some of them as 
canonical, or canonical in the homology class $[f_{g(v)}(M)]$, and  
then show that the zero section $f_{g(0)}(M)$ is a canonical cycle if, 
and only if, $s_g$ is constant.  
In \S4, we present a variety of examples of minimal unit, and  
minimal singular unit vector fields, some of canonical cycles and their 
scale deformations, 
and prove the characterization of
minimal unit vector fields on standard odd fimensional spheres indicated above,
which in dimension 3 we use to provide also a novel proof of the 
subject precursor \cite[THEOREM p.177]{glzi} of Gluck and Ziller, by employing
a Wirtinger like inequality between the volume forms of the 3d unit cycles 
involved.  
 
\section{Connections}
Our $M=M^m$s are connected closed smooth manifolds of dimension 
$m\geq 1$. Manifolds of tensor fields are provided with the compact open 
smooth topology, and subspaces with the induced subspace topology;  
$\mc{M}(M)$ is the manifold of Riemannian metrics of $M$, and for
$g \in \mc{M}(M)$, $\mc{M}_{[g]}(M)$ is the submanifold of metrics of $M$
in the conformal class of $g$.

\subsection{Connections on a vector bundle}
Let $E\stackrel{\pi}{\rightarrow} M$ be a vector bundle of rank $k$ with 
connection $\nabla^E$. 
In a local frame of sections $\{ e_1, \ldots, e_k\}$ over a coordinate 
neighborhood $U$ with coordinate functions $x=(x^1, \ldots, x^m)$, a section 
$s: M \rightarrow E$ can be written as 
$$
s= s^i e_i \, ,
$$
and $\{ \partial_{x^i},\partial_{s^j}\}$ yields a local frame for $TE$. 
Then we have that
$$
\nabla^E s = (\partial_{x^k} s^i + s^l \Gamma_{kl}^i ) e_i \otimes dx^k \, ,
$$
where the Christoffel symbols $\Gamma_{kl}^i$ are defined by 
\begin{equation} \label{eq1}
\nabla^E_{\partial_{x^k}}e_l= 
\Gamma_{kl}^i e_i \, .
\end{equation}
Thus 
$$
\nabla^E : C^{\infty}(E) \rightarrow C^{\infty}(E \otimes T^*M)\, .
$$

The connection $\nabla^E$ encodes a natural splitting
of $T_{p,s}(E)$ at each point $(p,s)\in TE$. The kernel of the  
smooth map $\pi_* :TE \rightarrow TM$ 
at any $(p,s(p))\in E$ defines the vertical subspace $V_{p,s(p)}E$ at the
point $(p,s(p))\in E_p=\pi^{-1}(p)\subset E$. Notice that $E_p$ is a vector 
space submanifold of $E$, so we have a natural identification 
$T_{p,s(p)}E_p \cong E_p$. On the other hand, if $x(t)$ is a curve on $M$ that 
starts at $p=x(0)$, the  curve $t\rightarrow s(x(t))$ is parallel along 
$x(t)$ if $s_{x(t)}\in T_{x(t)}E$ solves the parallel transport differential 
equation
$$ 
\nabla^E_{\dot{x}} s = 0 = \dot{x}^j( \partial_{x^j}s^i+ s^l \Gamma_{jl}^i)e_i
$$ 
with initial condition $s_{x(0)}$. We obtain a natural lift of 
$x(t)$ in $M$ to a curve 
$$ 
\gamma(t)=(x(t), s(x(t))) 
$$ 
in $E$, the velocity of which $\dot{\gamma}(t)\in T_{\gamma(t)}E$ 
is said to be a horizontal vector. We have that 
$$
\dot{\gamma}(t)=\gamma(t)_* \frac{\partial}{\partial t}\mid_{\gamma(t)}= 
\dot{x}^j \partial_{x^j}+ \dot{s}^i \partial_{s^i}
= \dot{x}^j \partial_{x^j}+ (-\Gamma_{jl}^i \dot{x}^j s^l) \partial_{s^i}
=\dot{x}^j (\partial_{x^j} -\Gamma_{jl}^i s^l \partial_{s^i})
=\dot{x}^j D_j  \, ,
$$
where
$$
D_j= \partial_{x^j} -\Gamma_{jl}^i s^l \partial_{s^i} 
\in T_{\gamma(t)}E
$$
is the covariant derivative operator on $T E$ in the direction $x^j$, which
is well-defined and independent of the curve $x(t)$, transversal to
$V_{x(t),s(x(t))}E$ for all $t$. Thus, 
$$
T_{p,s(p)}E= {\rm Vert}_{p,s(p)}E \oplus {\rm Hor}_{p,s(p)}E= \text{span$\{ 
\partial_{s^i}\} \oplus$ span$\{ D_j \}$\, .}
$$
Notice that $s\in E_p$ 
can be lifted to $T_{p,s}(E)$ either as a vertical vector tangent to the fiber 
$E_p$, or as a horizontal vector as above, induced by parallel translation. 
 
By using a a connection $\nabla$ on $T^{*}M$, or equivalently, on $TM$, and
$\nabla^E$, we produce a covariant differentiation 
$$
\nabla^{E,\nabla} : C^{\infty}(E\otimes T^*M) \rightarrow 
C^{\infty}(E \otimes T^*M\otimes T^*M)
$$
satisfying the Leibniz rule. 
If the Christoffel symbols of $\nabla$ are defined by 
$$
\nabla_{\partial_{x^i}}\partial_{x^j}= \gamma_{ij}^k \partial_{x^k} \, ,
$$
and so we have that
$$
\nabla_{\partial_{x^i}}dx^j= -\gamma_{ij}^k d{x^k} \, , 
$$
then, by (\ref{eq1}), we get 
that $\nabla^{E,\nabla}\nabla^E s =: \nabla \nabla s$ is given by
\begin{equation} \label{eq2}
\nabla \nabla s=   
\partial_{x^j}(\partial_{x^k} s^i + s^l \Gamma_{kl}^i ) e_i \otimes 
dx^j \otimes dx^k +
(\partial_{x^k} s^i + s^l \Gamma_{kl}^i )\Gamma_{ji}^p e_p \otimes dx^j
 \otimes dx^k - 
(\partial_{x^k} s^i + s^l \Gamma_{kl}^i )e_i \otimes \gamma_{jk}^l dx^l 
\otimes dx^j \, . 
\end{equation}
The curvature tensor $\Omega^{\nabla^E} \in \Lambda^2(M,{\rm End}\, E)$ 
of $\nabla^E$ is defined as
the commutator of covariant derivatives of a section, and the 
torsion tensor of $\nabla$ as  $T^{\nabla}=T^{\nabla \, k}_{ij}= 
\gamma_{ij}^k-\gamma_{ji}^k$.
By (\ref{eq2}), we obtain that
\begin{equation} \label{eq3}
\nabla_{\partial_{x^j}}\nabla_{\partial_{x^k}} s -
\nabla_{\partial_{x^k}}\nabla_{\partial_{x^j}} s = \Omega_{ljk}^{\nabla^E\, i} 
s^l e_i -
T_{jk}^{\nabla \, l} (\partial_{x^l} s^i + s^p \Gamma_{lp}^i ) e_i\, .  
\end{equation}
The curvature $\Omega^{\nabla_E}$ obstructs the integrability of the 
horizontal distribution (other than the zero section).

\begin{theorem} {\rm (Ambrose \& Singer \cite{amsi})} \label{le1} 
Let $\nabla^E$ and $\nabla$ be connections on $E$ and $TM$ as above,
with curvature and torsion $\Omega^{\nabla^E}$ and $T^\nabla$, respectively. 
Then:
\begin{enumerate}[label={\rm (\alph*)}]
\item If $s$ is the horizontal lift of $s_p \in E_p=\pi^{-1}(p)\subset E$, we 
have that
$$
\nabla_{\partial_{x^j}}\nabla_{\partial_{x^k}} s - \nabla_{\partial_{x^k}}
\nabla_{\partial_{x^j}} s = \Omega_{ljk}^{\nabla^E\, i} s^l e_i \, .
$$
\item If $T^\nabla \equiv 0$, then for any section $s$ of $E$ we have that  
$$
\nabla_{\partial_{x^j}}\nabla_{\partial_{x^k}} s - \nabla_{\partial_{x^k}}
\nabla_{\partial_{x^j}} s = \Omega_{ljk}^{\nabla^E\, i} s^l e_i \, .
$$
\end{enumerate}   
\end{theorem}

\subsection{The Sasaki metric on the tangent bundle of a Riemannian manifold}
If $g\in \mc{M}(M)$, its Levi-Civita connection $\nabla^{TM,g}$ provides 
a metric preserving torsion free connection on the rank $m$ bundle $E=TM$,
and an associated decomposition of $T(TM)$ into horizontal and vertical
subspaces at each point. We choose a normal frame
$\{ e_i \}$ on a neighborhood $U_p$ of $p$ in $M$, with coordinates
$x(q)=(x^1(q), \ldots , x^m(q))$,
and write locally a vector field $v$ on $M$ as 
$$
v(p)=\sum_i  v^i(p) e_{i \, p} \, .
$$
The pair $(x,v)=(x^1, \ldots, x^m, v^1, \ldots, v^m)$ defines local 
coordinates on $\tilde{U}_p:=\pi_M^{-1}(U_p)$ for $TM$, and a local 
chart for $T(TM)\stackrel{\pi_{TM}}{\rightarrow }TM$ on $\pi_{TM}^{-1}(
\tilde{U}_p)$ is defined by
$$
(x^1(p), \ldots, x^m(p), v^1(p), \ldots, v^m(p),
X^1, \ldots, X^m, V^1, \ldots, V^m) \, ,
$$
where
$$
s=\sum_{i}\left(X^i\partial_{x^i}\mid_{p,v(p)} +
V^i \partial_{v_i}\mid_{p,v(p)}\right) \in T_{(p,v(p))}TM 
$$
is a local section of $T(TM)$, so we have that
$$
\pi_{TM}(p,v(p),X^1, \ldots, X^m, V^1, \ldots, V^m) = 
(p, v(p)) \, , 
$$
and 
$$
(\pi_{TM})_{{*}_{p,v}} (p,v(p),X^1, \ldots, X^m, V^1, \ldots, V^m) = 
(p, X_{p}) \, ,
$$
respectively.
The vertical lift of 
$X_p\in T_pM$ to $T_{p,v}(TM)$ is given
by
$$
X^{ver}_{p,v(p)}=X^1 (p)\partial_{v^1}\mid_{v(p)} +\cdots
+ X^m (p) \partial_{v^d}\mid_{v(p)}\, , 
$$
while if we take a curve in $M$ that starts at $p$ with initial 
velocity $X_p$, and by using the connection (\ref{eq1}) of $g$, 
the parallel transport of $X_p$ along the curve
defines its horizontal lift 
$$
X^{hor}_{p,v(p)}=(X^1\partial_{x^1} -X^j \Gamma_{jl}^1 v^l \partial_{v^1})+
\cdots + (X^m \partial_{x^d} -X^j \Gamma_{jl}^n v^l \partial_{v^m})  \, ,  
$$
and we have that $(\pi_{TM})_{* p,v(p)}X^{hor}_{p,v(p)}=X_p$, so a vector field 
$X$ on $M$ lifts in this manner to a unique horizontal vector field $X^{hor}$ 
in $TM$ that is $\pi_{TM}$ related to it. 
We have that
$$ 
T_{p,v(p)}(TM)= {\rm Vert}_{p,v(p)}(TM) \oplus 
{\rm Hor}_{p,v(p)}(TM)= \text{span$\{ 
\partial_{v^i}\} \oplus$ span$\{ \partial_{x^j}-\Gamma^i_{jl}v^l \partial_{v^i}
\}$\, ,}
$$
and the coordinate decomposition of
$\Xi=
(X^1, \ldots, X^m, V^1, \ldots, V^m)_{p,v(p)} \in T_{p,v(p)}TM$ into horizontal
and vertical components is given by   
$$
(X^1, \ldots, X^m, -X^j \Gamma_{jl}^1 v^l , \ldots, -X^j \Gamma_{jl}^m v^l )_{
p,v(p)} + ( 0, \ldots, 0, V^1+X^j\Gamma_{jl}^1 v^l , \ldots,
V^m+X^j \Gamma_{jl}^m v^l )_{p,v(p)} \, . 
$$
Notice that the connection operator $K$ of $\nabla^{TM,g}$ identifies the
vertical component of $\Xi$ with an element of the vector space $T_pM$ as
a submanifold of $TM$, and $T_{p,v(p)}T_pM \cong T_pM \hookrightarrow 
T_{p,v(p)}TM$ equals the kernel of $(\pi_{TM}){*_{p,v(p)}}$. Indeed, if 
$$
T_{p,v(p)}T_pM=\{ (p,v(p), 0, V): \; V\in \mb{R}^{m}\} \, ,
$$ 
then the mapping 
$$ 
\iota_{p,v(p)}:  T_{p,v(p)}T_pM \ni (p,v(p),0,V) \rightarrow
(p,V) \in T_p M 
$$ 
provides the indicated isomorphism $K(\Xi)=\iota_{p,v(p)}(\Xi^{ver})$ 
with $T_pM$.

The Sasaki metric $g_S$ on $TM$ is defined by declaring 
$\{e_1^{hor}, \ldots, e_m^{hor},e_1^{ver}, \ldots, e_m^{ver}\}$ to be a
normal frame field on $\pi^{-1}(U_p)$. Its dual coframe is given by 
$$
\pi_M^* e^1, \ldots, \pi_M^* e^m, Dv^1, \ldots, Dv^m\, ,
$$ 
where $e^i(e_j)=\delta^i_j$ and $Dv^i= dv^i + v^l\pi_M^* (\Gamma_{jl}^ie^j)$.
Thus, as in its original definition \cite{sas}, $g_S$ can be equally expressed
as 
$$
g_S(\Xi,\Xi') =  g((\pi_{TM}){*_{p,v(p)}}(\Xi),(\pi_{TM}){*_{p,v(p)}}(\Xi'))
+g(K(\Xi),K(\Xi')) \, , 
$$
where $\Xi, \Xi' \in T_{p,v(p)}(TM)$, or as the $(0,2)$ symmetric tensor
$$
g_S = (\pi_M^* e^i)^2 + (Dv^i)^2 \, .  
$$
It is uniquely characterized as the metric that makes the fibration
$$
\pi : (TM,g_S) \rightarrow (M,g)
$$
a Riemannian submersion, is the constant metric $g_p$ along each tangent 
space fiber $T(T_pM) \cong  T_p M \hookrightarrow T(TM)$, and along the 
horizontal space defined by the Levi-Civita connection of $g$, it is the metric 
coming from the base by projection. Any fiber of this bundle is a totally 
geodesic flat submanifold of the total space.
The unit sphere bundle $\iota_{S^1(TM)}: S^1(TM)\hookrightarrow T(TM)$ 
inherits the metric $\iota_{S^1(TM)}^* g_S$, which for notational convenience
we refer to as $g_S$ also. The fibration
$$ 
(S^1(TM),g_S) \rightarrow (M,g) 
$$ 
is a Riemannian submersion with totally geodesic fibers also.
Expressions for the Levi-Civita connection and curvature of $g_S$, 
of relevance to us here, can be found in the literature, for instance,
\cite{mutr,kowa,guka}.
 
The flow map $\varphi_t : TM \rightarrow TM$ of the 
tautological vector field 
$$
V\mid_{p,v(p)}=\left( v^1\partial_{x^1}
-v^j \Gamma_{jl}^1 v^l \partial_{v^1}\right) +\cdots
+ \left( v^m \partial_{x^m} 
-v^j \Gamma_{jl}^m v^l \partial_{v^m} \right)\mid_{p,v(p)} 
$$
is defined by
$\varphi_t(p,v(p))=(\gamma(t),\dot{\gamma}(t)) \in TM$,  
where $\gamma(t)$ is the unique geodesic through $p$ with initial velocity 
$\dot{\gamma}(t)\mid_{t=0}=v(p)$. 
Thus, $\varphi_t(S^1(TM))
\subset S^1(TM)$ and $V\mid_{S^1(TM)}$ is tangent to $S^1(TM)$.   

The oriented manifolds $TM$ and $T^*M$ are naturally identified with each
other by the metric $g$. The projectivization 
$(T^*M\setminus \{0_{TM}\})/\mb{R}_{+}$ identifies
with $S^1(T^* M)$, and by metric duality with $S^1(TM)$. If 
$\lambda$ and $\omega=d\, \lambda$ are the 
canonical one 
and symplectic forms of  $T^*M$, respectively, 
the 
duality makes 
$\omega$  corresponds to a symplectic form $\omega^{\sharp}$ on $TM$.
The restriction of $\lambda$ defines a contact form $\sigma$ on 
$S^1(T^*M)$ whose kernel is the contact bundle $\mc{D}_{\sigma}$.
The Reeb vector field 
$v^*_{\sigma} \in C^{\infty}(S^1(T^*M), T(S^1(T^*M)))$ 
of the contact structure is defined uniquely by the conditions
$$
\sigma(v^*_{\sigma})=1\, , \quad d\sigma \hok v^*_{\sigma}=0 \, ,
$$
and the metric identifies it with the restriction to $S^1(TM)$ of the 
tautological field $v_{\sigma}=V_{p,v_{\sigma}(p)}$ that generates the 
geodesic flow.

Under a conformal change $e^{2u}g$ of $g$, $g_S$ changes conformally to 
$e^{2u \circ \pi_M}g_S$, and so the Levi-Civita connection 
$\nabla^{TM,g_S}$ of $g_S$ is invariant under rescaling by $l$ of $g$. Under 
this
rescaling, the Riemannian connection $\nabla^{TM,g}$ changes to a connection
$\nabla^{TM,l,g}$ with a well-defined adiabatic limit 
$\nabla^{TM,0,g}$ as $l\rightarrow \infty$ that collapses the fibers.
We have that  $\nabla^{TM,g_S}$ and $\nabla^{TM,0,g}$ are CS equivalent
\cite[Theorem 3.2]{sisu}, so if $P$ is any invariant polynomial of degree $l$, 
$l\leq m$, we not only have equality of the characteristic classes 
$[P(\Omega^{\nabla^{TM,0,g}})]$ and $[P(\Omega^{\nabla^{TM,g_S}})]$, but 
equality of the characteristic forms $P(\Omega^{\nabla^{TM,0,g}})$ and 
$P(\Omega^{\nabla^{TM,g_S}})$ themselves.

\section{Vector field cycles in $TM$} 
The geometry of each $f_{g(v)}(M)$ in (\ref{egv}) is encoded in the 
Whitney sum decomposition 
$$
\iota_v^* T(TM) = T(f_{g(v)}(M))\oplus \nu(f_{g(v)}(M))\, , 
$$
where $\iota_v : f_{g(v)}(M) \hookrightarrow TM$ is the inclusion map. 
For if we let $\alpha^{TM,g_S}_{f_{g(v)}}$ and $H^{TM,g_S}_{f_{g(v)}}$ be the 
second fundamental form and mean curvature vector of the embedding 
$f_{g(v)}$, the scalar curvature $s_{f_{g(v)}}$ of the metric $g(v)$ 
satisfies the relation \cite{gracie}
\begin{equation} \label{sca1}
s_{f_{g(v)}}= \sum_{i,j} K^{TM,g_S}(e^{f_{g(v)}}_i,e^{f_{g(v)}}_j) + 
\| H^{TM,g_S}_{f_{g(v)}}\|^2-\| 
\alpha^{TM,g_S}_{f_{g(v)}}
\|^2 \, , 
\end{equation}
where the first summand on the right is the exterior scalar curvature of
$g_S$ on $f_{g(v)}(M)$. 
We analyze how the volume of $g(v)$ is encoded in the 
vector field $v$, and position ourselves towards the study of conformal
deformations of $f_{g(v)}(M)\hookrightarrow TM$.

\subsection{Minimal vector fields}
If $d\mu_{f_{g(v)}}$ is the Riemannian volume measure of $g(v)$, the family 
(\ref{egv}) has an associated family of volume functionals given by 
\begin{equation} \label{vgv}
C^{\infty}(M;TM)\ni v \mapsto 
\mu^{TM}_{g(v)}(v)=\int_{f_g(v)(M)} d\mu_{f_{g(v)}} \, . 
\end{equation}
Along a path $t\rightarrow v_t$ of deformations of $v$, the form 
$d\mu_{f_{g(v_t)}}$ of $f_{g(v_t)}(M)$ varies as 
\begin{equation} \label{veqv}
\frac{d}{dt}d\mu_{f_{g(v_t)}}=\left( {\rm div}(T^{\tau}) - \< T^{\nu},
H^{TM,g_S}_{f_{g(v_t)}}\> \right) d\mu_{f_{g(v_t)}}\, ,
\end{equation}
where $T=T^{\tau}+T^{\nu}$ is the decomposition of the variational vector
field $T=f_{g(v_t)\, *} \partial_t$ into tangential and normal components to 
the embedded submanifold in $(TM,g_S)$.  Thus, the critical points of 
(\ref{vgv}) are vector fields $v$ such that the embedding $f_{g(v)}(M)$ is 
minimal. Such a $v$ is by definition a {\it minimal vector field} of $(M,g)$. 
We pursue their geometric characterization, and show that with this degree of 
generality, they are scarce, a fact that in on itself provides information on 
the geometry of $(M,g)$. 
\smallskip

The symmetric and antisymmetric decomposition of the $(0,2)$-tensor 
$C_v:=\nabla^g\, v$ is given by
\begin{equation} \label{iso}
C_v(X,Y)=\frac{1}{2}\mc{L}_vg(X,Y)+
\omega_v(X,Y)\, ,
\end{equation}
where $\mc{L}_v g$ is the Lie derivative of the metric tensor $g$ in the 
direction of $v$, and if $\tilde{v}=g(v,\, \cdot \,)$ is the one form metric 
dual to $v$, $\omega_v:=d\tilde{v}$ is the vorticity of $v$. The vector 
field $v$ is said to be conformal if its flow generates a one 
parameter group of conformal transformation, and so we have that 
$\frac{1}{2}\mc{L}_v g = u_v g$ for some scalar function $u_v$. If we view 
$C_v$ as the section of ${\rm End}\, T M$ given by the family of   
$(1,1)$ tensors $X \rightarrow \nabla_X^g v$ on the 
tangent space, the   
family of nonnegative operators $C_v^t C_v$ admits locally an 
orthonormal frame of eigenvectors, any one of which induces in turn 
a local orthonormal tangent frame of $f_{g(v)}(M)$. 
  
In what follows, $R^g$ will stand for the Riemannian curvature tensor
$R^g(X,Y)Z=(\nabla^g_X \nabla^g_Y- \nabla^g_Y \nabla^g_X- \nabla^g_{[X,Y]})Z$
of the metric $g$. 

\begin{lemma} \label{le2}
Let $\{ e_i\}$ be any local orthonormal frame of eigenvectors
of $C_v^t C_v$ about $p\in M$, 
and let $\{ e_i^{f_{g(v)}}:=c_i(e_i^{hor}+(\nabla^g_{e_i}v)^{ver})\}_i$,
$c_i^{-2}= 1+\| \nabla^g_{e_i}v\|^2=:1+\lambda^2_i(v)$, 
be the associated orthonormal tangent frame of $f_{g(v)}(M) \hookrightarrow 
(TM,g_S)$. Then 
\begin{equation} \label{mcv}
H_{f_{g(v)}}^{TM,g_S}= \left( \sum_{i}  
c_i^2\left((\nabla^g_{e_i} e_i)^{hor}+ (\nabla^g_{e_i}(\nabla^g_{e_i} v))^{ver}+
(R^g(v,\nabla^g_{e_i} v)e_i)^{hor}\right)\right)^{\nu=\nu_{TM}}\, ,
\end{equation}
where $X^\nu$ is the component of $X \in T(TM)\mid_{f_{g(v)}(M)}$ 
normal to $f_{g(v)}(M)$.  
In particular, $v$ is minimal if, and only if, we have that
$$
\sum_{i} c_i^2\left(\nabla^g_{e_i}\nabla^g_{e_i} -
\nabla^g_{\nabla^g_{e_i} e_i} -\nabla^{g}_{R^g(v,\nabla^g_{e_i} v)e_i}\right)v 
=0 \, ,
$$
and $v$ satisfies this equation if, and only if, $v$ is either parallel or 
identically zero.    
\end{lemma}

{\it Proof}. If $X$ and $Y$ are local vector fields in $M$, the vector fields 
$\tilde{X}=X^{hor}+(\nabla^g_X v)^{ver}$ and $\tilde{Y}=Y^{hor}+
(\nabla^g_Y v)^{ver}$ in $TM$ are tangent to $f_{g(v)}(M)$, and we have 
that
$$
\nabla^{g_S}_{\tilde{X}}\tilde{Y}=
(\nabla^g_X Y)^{hor}-\frac{1}{2}R^g(X,Y)v)^{ver}+
(\nabla^g_X (\nabla^g_Y v))^{ver}+\frac{1}{2}R^g(v,\nabla^g_Y v)X)^{hor}+
\frac{1}{2}R^g(v,\nabla^g_X v)Y)^{hor}\, .
$$
Thus, in the orthonormal frame $\{ e_i^{f_{g(v)}}\}$ of 
$(f_{g(v)}(M),g(v))$, we have that  
$$
\nabla^{g_S}_{e_i^{f_{g(v)}}}e_i^{f_{g(v)}}=c_i^2\left( 
(\nabla^g_{e_i} e_i)^{hor}+ (\nabla^g_{e_i} (\nabla^g_{e_i} v))^{ver}+
(R^g(v,\nabla^g_{e_i} v)e_i)^{hor} - 
 e_i(\lambda_i) \left( e_i^{hor}+(\nabla^g_{e_i}v)^{ver} \right)
\right) \, ,
$$
with the last term on the right above tangential.
The first statement follows.

At the point $(p,v(p))$ of $f_{g(v)}(M)$, we choose the set of normal
vectors $\{ N_i = -((\nabla^{g}v)^t(e_i))^{hor}+ e_i^{ver}\}_i$.  Although not 
necessarily orthogonal, $\{ N_i\}_i$ forms a basis of 
$\nu_{(p,v(p))}(f_{g(v)}(M))$.
The projection of $H_{f_{g(v)}}^{TM,g_S}$ along $N_j$ has magnitude 
$$
\< \sum_i c_i^2\left(\nabla^g_{e_i}\nabla^g_{e_i} -
\nabla^g_{\nabla^g_{e_i} e_i} -\nabla^{g}_{R^g(v,\nabla^g_{e_i} v)e_i}\right)v 
,e_j\>\| N_j\|^{-1} \, , 
$$
hence, $v$ is minimal if, and only if, 
$$
\sum_i c_i^2\left(\nabla^g_{e_i}\nabla^g_{e_i} - \nabla^g_{\nabla^g_{e_i} e_i}
 -\nabla^{g}_{R^g(v,\nabla^g_{e_i} v)e_i}\right)v = 0 \, ,  
$$
as stated.

Since $C_v^t C_v (e_i)=\lambda_i^2(v) e_i$ where $\lambda_i^2(v)= 
\| \nabla^g_{e_i}v\|^2$, we have that
$$
d\mu_{f_{g(v)}^{*}g(v)}= (\prod_i (1+\lambda_i^2(v) )^{\frac{1}{2}}) d\mu_g \, ,
$$
and under the scaling $t \rightarrow tv$, $\lambda_i^2(tv) \rightarrow t^2 
\lambda_i^2(v)$. If at least one of the $\lambda_i$s is nonzero anywhere,   
$$
\frac{d}{dt}\int d\mu_{f_{g}(tv)^{*}g(tv)}\mid_{t=1}= 
\int \left( \frac{d}{dt}(\prod_i (1+\lambda_i^2(tv))^{\frac{1}{2}}) 
\right) d\mu_g\mid_{t=1}\neq 0\, .
$$

If $C_v^t C_v=(\nabla^g v)^t \nabla^g v=0$ then $v$ is a minimal vector
field. Thus, if $(M,g)$ carries a nontrivial parallel vector field $v$, then 
for any
$t \in \mb{R}$, we have that $(1-t)v$ is minimal and $\mu_{f_g((1-t)v)}(M)=
\mu_g(M)$ is the absolute minimizer of (\ref{vgv}). 
Generically, $(M,g)$ does not carry parallel vector fields, the 
absolute minimizer of (\ref{vgv}) is rigid, and the minimum is achieved 
only by the embedding of $M$ as the zero section of $TM$.
\qed

\subsection{Minimal unit vector fields} 
\label{ss22}
Suppose now that $v\in C^{\infty}(M;TM\setminus 0_{TM})$. Then 
$v_g = v/\| v\|_g$ is a unit vector field on $M$, the geodesible flow of $v$ 
is given by the integral curves of $v_g$, and 
$f_{g(v_g)}(M) \hookrightarrow S^1(TM)$. In this context, it is pertinent to
consider the family of embeddings (\ref{egvU}), and the associated functional
\begin{equation} \label{eq16}
C^{\infty}(M;S^1(TM))\ni v \; {\displaystyle \mapsto \mu^{S^1(TM)}_{g(v)}(v)= 
\int_{f_{g(v)}(M)} d\mu_{f_{g(v)} (M)}} \, ,  
\end{equation}
a critical point $v$ of which is said to be a {\it minimal unit vector 
field}. 
Along variations $v_t$ of $v$ in $S^1(TM)$ we now have that
$$ 
\frac{d}{dt}d\mu_{f_{g(v_t)}}=\left( {\rm div}(T^{\tau}) - \< T^{\nu},
H^{S^1(TM),g_S}_{f_{g(v_t)}}\> \right) d\mu_{f_{g(v_t)}}\, ,
$$
so $v$ is a minimal unit vector field if, and only if, 
$H^{S^1(TM),g_S}_{f_{g(v)}}$ vanishes.

When $\| v\|_g=1$, any $e$ orthonormal to $v$ is an eigenvector of
$C_v^tC_v$ of eigenvalue $\lambda^2=\|\nabla^g_e v\|^2$ and $C_v(e)=\lambda 
e'$ is orthogonal to $v$, with $e'=e$ if $\lambda=0$. 

\begin{theorem} \label{math}
If $v$ is a unit vector field in $(M^{1+m},g)$, let 
$\{ e_0=v, \{ e_i\}_{i=1}^m\}$ be any local orthonormal frame about
$p\in M$ {\rm (}the set $\{ e_i\}_{i=1}^m$ is then a local orthonormal family 
of eigenvectors of $C_v^t C_v${\rm )}, and let 
$\{  e_i^{f_{g(v)}}:=c_i(e_i^{hor}+ (\nabla^g_{e_i}v)^{ver}_i\}$, 
$c_{i}^{-2}= 1+\| \nabla^g_{e_i}v\|^2 = 1+\lambda^2_i$,  
be the associated orthonormal tangent frame of $f_{g(v)}(M) \hookrightarrow 
(S^1(TM),g_S)$. Then $v$ is a minimal unit vector field if, and only if, 
\begin{equation} \label{eq10}
\sum_{i\geq 0} c_i^2\left(\nabla^g_{e_i}\nabla^g_{e_i} -
\nabla^g_{\nabla^g_{e_i} e_i} -\nabla^{g}_{R^g(v,\nabla^g_{e_i} v)e_i}\right)v 
=\lambda v \, ,   
\end{equation}
where $\lambda = - \sum_{i\geq 0} c_i^2 \lambda_i^2=-\sum_{i\geq 0} 
\lambda_i^2/(1+\lambda_i^2)$. A vector field $v$ that is both, minimal and 
minimal unit, must be a parallel vector field of unit norm.  
\end{theorem}

{\it Proof}.  If $X$ is a vector field in $M$, the vector field 
$\tilde{X}=X^{hor}+(\nabla^g_X v)^{ver}$ is tangent to $f_{g(v)}(M)$, and
since $\| v\|_g^2=1$, it is tangent to $S^1(TM)$ also.  We obtain that
$$
\nabla^{S^1(TM),g_S}_{\tilde{X}}\tilde{X}= (\nabla^g_X X)^{hor} + 
(\nabla^g_X (\nabla^g_X v))^{ver}- \<  \nabla^g_X (\nabla^g_X v),v\> v^{ver} 
+R^g(v,\nabla^g_X v)X)^{hor}\, , 
$$
and therefore, 
$$
\begin{array}{rcl}
\nabla^{S^1(TM),g_S}_{e_i^{f_{g(v)}}}e_i^{f_{g(v)}} & = & c_i^2\left( 
(\nabla^g_{e_i} e_i)^{hor}+ (\nabla^g_{e_i} (\nabla^g_{e_i} v))^{ver}
+ (R^g(v,\nabla^g_{e_i} v)e_i)^{hor} - 
e_i(\lambda_i) \left( e_i^{hor}+(\nabla^g_{e_i}v)^{ver} \right) \right)\\
& & -\<  \nabla^g_{c_i e_i} (\nabla^g_{c_i e_i} v),v\> v^{ver} \, .
\end{array}
$$
Hence, 
$$
H_{f_{g(v)}}^{S^1(TM),g_S} = \left( \sum_i c_i^2\left( 
(\nabla^g_{e_i} e_i)^{hor}+ (\nabla^g_{e_i} (\nabla^g_{e_i} v))^{ver}
+ (R^g(v,\nabla^g_{e_i} v)e_i)^{hor} \right)\right)^{\nu}\, ,  
$$
where $X^\nu$ is the component of $X$ normal to $v$ and
$f_{g(v)}(M)$ at $(p,v(p))$.

Since $\< C_v(e_0),e_0\>  = 0$, the 
set $\{ e_0^{ver}, \{ N^{TM}_i=-((\nabla^{g}v)^t(e_i))^{hor}+
e_i^{ver}\}_{i\geq 1}\}$ is a basis of the normal bundle of the submanifold 
$f_{g(v)}(M)$ in $TM$ at each $(p,v(p))\in S^1(TM)$ where defined, 
and the set
$\{ N^{S^1(TM)}_i=-((\nabla^{g}v)^t(e_i))^{hor}+e_i^{ver}\}_{i\geq 1}$ 
is a basis of its normal bundle in $S^1(TM)$. The vector 
$H_{f_{g(v)}}^{S^1(TM),g_S}$ is null if, and only if, the projection 
of the vector in the expression above along any of the $N^{S^1(TM)}_i$s is 
the null vector in $T_p M$, and so the said vector 
is proportional to $v^{ver}$. The desired result follows, 
the explicit form of the scalar $\lambda$ being a simple consequence of the
unitarity of $v$. Notice that $\lambda^2= \| H_{f_{g(v)}}^{TM,g_S}\|^2$, 
and so the minimal unit vector field $v$ is minimal if, and only if,   
$\nabla^g_{e_j}v=0$, $0\leq j\leq m$, and so $v$ must then be
parallel.  
\qed

The Euler characteristic of the manifold obstructs the existence of a 
nonsingular vector field on it, but the tangent bundle of any Riemannian 
manifold always has a nonzero section over the complement of a finite number 
of points, and the Hopf index theorem says that the Euler number of the 
bundle is the sum of the local degrees of this section at the said singular 
points. Since the condition encoded in (\ref{eq10}) is local, we can extend the
concept discussed above for vector fields with singularities, zeroes or 
otherwise, in a closed subset of the manifold: We say 
that a vector field $v$ {\it is a minimal unit singular vector field} if its
singular set $Z(v)$ is a closed subset of $M$, and on $M\setminus Z(v)$,
identity (\ref{eq10}) holds relative to any frame $\{ e_0=v/\|v\|_g,  
\{ e_j\}_{j>0}\}$ where $\{ e_j\}_{j>0}$ is a family of eigenvectors of 
$C^t_{e_0}C_{e_0}$.   
When $Z(v)$ is not empty, $f_g(v/\|v\|_g)(M\setminus Z(v))$ is no 
longer a closed submanifold of $(TM,g_S)$, though under some hypothesis on 
$Z(v)$, $f_g(v/\|v\|_g)(M\setminus Z(v))$ can be interpreted as an 
$m$ current on $TM$ that encodes (some) topological properties of $M$. 
The parallel fields with one singularity of \cite{pede} and \cite{bogi2}
(push-forward of vector fields on the blow-up of $M$ at the singular point),  
both constructed to produce convergent volume minimizing sequences 
of unit vector fields, are natural examples of this situation. 


\begin{theorem} \label{th4}
If $v$ is a minimal unit vector field on $(M^{1+m},g)$, 
the leaves of the foliation of $M$ it defines are geodesics, and its flow is 
by isometries.
\end{theorem}

{\it Proof}. If $m=0$, or $v$ is parallel, there is nothing to prove, so we 
assume otherwise. Thus, if $\{ e_0=v, e_1, \ldots, e_m\}$ is a local  
orthonormal frame, the set $\{e_i\}_{i\geq 1}$ is an orthonormal family of 
eigenvectors of $C_v^tC_v$, and identity (\ref{eq10}) holds with $\lambda < 0$. 

Suppose that $\nabla^g_v v \not \equiv 0$. Since at any point $p\in M$ where 
$\nabla^g_v v \neq 0$ the direction of $\nabla^g_v v$ is uniquely defined, 
there exists a unit vector field $e_r$, defined on some maximal open set
$\mc{U}\ni p$, such that
\begin{equation} \label{to1} 
\nabla^g_v v =c \, e_r\, ,  
\end{equation}
for a positive scalar function $c$.

We take any completion of the pair $\{ v, e_r\}$ to a local 
orthonormal frame $\{ e_0:=v, e_1:=e_r,e_2, \ldots, e_m\}$. Since $\| 
\nabla_v^gv\|^2 \neq 0$ and the dimension of the image space of $C_v$ is at 
most $m$, in the family $\{ e_i\}_{i\geq 1}$ of eigenvectors of $C_v^tC_v$ 
there is at least one of eigenvalue $0$.  
If we write    
\begin{equation} \label{neq}
\nabla^g_v \nabla^g_v v = v(c)e_r +c\nabla^g_v e_r = 
v(c)e_r-c^2 v +c\sum_{j\geq 2} \beta_j e_j \, ,
\end{equation}
where $\beta_j e_j= \< \nabla^g_v e_r ,e_j\> e_j$, after some minor 
manipulations,
(\ref{eq10}) can be written equivalently as
$$
\begin{array}{lcl}
{\displaystyle \frac{1}{1+c^2}(v(c) e_r + c\sum_{j\geq 2} \beta_j e_j -
c( \nabla_{e_r+R^g(v,e_r)v})v)}+ & &  \\
\sum_{j\geq 1}{\displaystyle \frac{1}{1+\lambda_j^2}(\nabla^g_{e_j}
\nabla^g_{e_j}v-
\<\nabla^g_{e_j}\nabla^g_{e_j}v,v\> v - \nabla_{\nabla_{e_j}e_j+R^g(v,
\nabla_{e_j}v)e_j}v )} &  = & 0 \, , 
\end{array}
$$  
and since at any fixed $q\in \mc{U}$ we can assume that 
$\{e_r, e_2, \ldots, e_m\}$ is such that 
$$
\nabla^g_{e_r} e_r \mid_{q}=
\nabla^g_{e_r} e_k \mid_{q}=0=\nabla_{e_{k'}}e_k\mid_q \, , \quad k,k'\geq 2\, ,
$$
the $e_{k\geq 2}$th component of the vector on the left side of the identity 
above can be seen to be $\beta_ke_k\mid_q$, which then must be zero, and by
the arbitrariness of $q$, we thus conclude that $\nabla^g_v e_r = -c v$, and 
then by (\ref{neq}) that $v(c)=0$, while $\nabla^g_{e_r}v\mid_q=
\left(\sum_{j\geq 1}\<\nabla^g_{e_r}v,e_j\>e_j\right)\mid_q=0$, and 
$\nabla^g_{e_r}e_r\mid_q= \left(\< \nabla^g_{e_r}e_r,v\>v+ \sum_{k\geq 2}
 \< \nabla^g_{e_r}e_r,e_k\> e_k\right)\mid_q=0$, so
$\nabla^g_{e_r}v=0=\nabla^g_{e_r}e_r$. 

Since $\nabla^g$ is torsion free, we then have that   
$$
\nabla^g_v e_r -\nabla^g_{e_r}v=-c\, v= [v,e_r] \, ,
$$
and the distribution of planes spanned by $\{v, e_r\}$ is integrable. 
We let $\Sigma$ be the unique integral manifold of this distribution through 
$p$, and provide it with the metric $g_{\Sigma}=i^* g$ induced by the 
metric $g$ of the  ambient space. Then $(\Sigma,g_{\Sigma})$ is 
a totally geodesic flat surface embedded into $(M,g)$. 

Restricted to $\Sigma$, $v$ defines a minimal unit vector field.   
In the orthonormal frame $\{ v, e_r\}$, equation (\ref{eq10}) reduces to the 
identity
\begin{equation} \label{to2}
\nabla^g_v \nabla^g_v v = - c^2 v  \, ,   
\end{equation}
(which is the restriction to $\Sigma$ of the equivalent equation  
$v(c)e_r -c^2 v = - c^2 v$ that holds on 
the entire $\mc{U} \supset \Sigma$).
Since   
$v(c)=0$, $c$ is a positive constant along integral curves of $v$, 
which are therefore entirely contained in $\Sigma \subset \mc{U}$. 
On the other hand, if we apply $\nabla^g_{e_r}$ to both sides of (\ref{to2}),
by the flatness of the metric, we obtain that
$$
\begin{array}{rcl}
\nabla^g_{e_r}\nabla^g_v \nabla^g_v v  & = & - e_r(c^2) v  \\ 
(\nabla^g_{v}\nabla^g_{e_r}+\nabla^g_{[e_r,v]}) \nabla^g_v v  & = &  \\ 
-2c^3 v & = & 
\end{array} \, , 
$$
and so $e_r(c)=c^2$ (by applying $\nabla^g_{e_r}$ in  
(\ref{to1}), we would obtain similarly that $c^2 e_r = e_r(c)e_r$, and arrive 
at the same equation for $e_r(c)$ over the entire set $\mc{U}$). Thus, along 
the integral curves of $e_r$, $c$ stays positive, hence the entire integral 
curve is contained in $\Sigma \subset \mc{U}$ also, and blows up to 
$+\infty$ in finite 
time. This contradicts the fact that $c^2$ is the restriction to $\Sigma$ of 
the globally defined smooth function $\| \nabla_v^gv \|_g^2$. 
Therefore, $\nabla^g_v v = 0$, the flow of $v$ is by geodesics, and the 
geodesible flow of $v$ is given by its integral curves.       

Along the orbits of the geodesible flow, the orthogonal plane field is 
invariant \cite[Theorem ii, \S1]{sull}. If we work on 
a tubular neighborhood of a geodesic segment of the foliation, and $V$ is a
vector field tangent to the foliation, and $Y$ is (a local) unit orthonormal, 
we have that
$$
\begin{array}{rcl}
0=\frac{1}{2}Vg(Y,Y) & = & g(\nabla_V^gY,Y)= g(\nabla_Y^gV+[V,Y],Y) \\
& = & Yg(V,Y)-g(V,\nabla_Y^gY)+g([V,Y],Y) \\
 & = & g([V,Y],Y) \, ,
\end{array}
$$
and the metric $g$ is bundle like relative to the foliation.  If $X$ is any 
vector field in the invariant orthogonal plane field, we then have that 
$d\tilde{v}(v,X)=0$  \cite[Theorem i,\S1]{sull}, and by (\ref{iso}) it follows
that $\mc{L}_vg(v,X)=g(\nabla_v v,X)+g(\nabla_X v, v) = 0$. Since the
flow is by geodesics, by (\ref{iso}) also, we have that $g(\nabla_vv,v)=0=
\mc{L}_vg(v,v)$ and that for any pair $X,Y$ of vector 
fields in the invariant orthogonal plane, 
$\mc{L}_vg(X,Y)=g(\nabla^g_Xv,Y)+g(\nabla^g_Yv,X)= 
-g(v,\nabla^g_XY+\nabla^g_YX)=0$. Thus, $\mc{L}_v g=0$, $v$ is a 
Killing field, and the flow is by isometries.  
\qed

\begin{remark}
No closed Riemannian manifold $(M,g)$ with metric $g$ of negative 
Ricci curvature carries minimal unit vector fields, as a Killing field is then
trivial \cite{boch}. This is consistent  with the fact that an 
$(M,g)$ of negative sectional curvature cannot carry this type of vector 
fields either, as such a manifold does not admit Riemannian flows at 
all for, if one such were to exist, the Gromov simplicial volume invariant of
$M$ would be zero \cite{grom}, \cite[Corollaire 5]{carr2}. In these manifolds,
any minimal singular unit vector field must have nontrivial singular set, and 
the limit of convergent volume minimizing sequences of unit vector fields 
must have singularities.   
\end{remark}
\smallskip

\subsection{Canonical vector field cycles}
Each $f_{g(v)}(M)$ is a closed submanifold of $TM$, and so if $M$ itself is 
oriented, its Poincar\'e dual and the Thom class of its normal bundle can be 
represented by the same form in $H^{*}_{c,ver}(TM)\cong H^{*-m}(M)$, in 
particular, by the Poincar\'e dual of the zero section of the tangent bundle; 
if $M$ is nonorientable, the same remains true by a form in
$H^{*}_{c,vert}(TM)\cong H^{*-m}(M;\Lambda)$, $\Lambda$ the bundle of 
1-densities over $M$. In the set of all of these singular homology classes, 
we seek ``canonical cycles'' exploiting the conformal geometry of the metric 
$g_S$ on $TM$, and the relation of these cycles to the vector fields on $M$ 
themselves, and their scale deformations. The deformations are of particular 
significance as recorders of the fact that the compactly supported 
cohomology of $TM$ is not invariant under homotopy equivalence, but it is so 
under diffeomorphism.

By (\ref{sca1}), the total scalar curvature functional decomposes as  
\begin{equation} \label{eq11}
\begin{array}{rcl} 
{\displaystyle S_{f_{g(v)}}(M):= \int s_{f_{g(v)}} d\mu_{f_{g(v)}} } 
& =  &  
{\displaystyle \Theta_{f_{g(v)}}(M):= \int_{f_{g(v)}(M)} \sum K^{TM,g_S}(
e^{f_{g(v)}}_i,e^{f_{g(v)}}_j) d\mu_{f_{g(v)}}} \vspace{1mm} \\
& & +{\displaystyle \Psi_{f_{g(v)}}(M) := \int_{f_g(v)(M)} \| H^{TM,g_S}_{
f_{g(v)}}\|^2 d\mu_{f_{g(v)}} }  \vspace{1mm} \\ & & 
-{\displaystyle \Pi_{f_{g(v)}}(M):= \int_{f_{g(v)}(M)} \| \alpha^{TM,g_S}_{
f_{g(v)}}\|^2 d\mu_{f_{g(v)}} } \, , \vspace{1mm} \\
& = & {\displaystyle  \mc{W}_{f_{g(v)}}(M) - \mc{D}_{f_{g(v)}}(M)}  \, ,
\end{array}
\end{equation}
where if $m\geq 2$, 
\begin{equation} \label{eq12}
\begin{array}{rcl}
\mc{W}_{f_{g(v)}}(M^m) & := & {\displaystyle 
\frac{m}{m-1}\Theta_{f_{g(v)}}(M)+\Psi_{f_{g(v)}}(M) }\, ,\vspace{1mm}\\
\mc{D}_{f_{g(v)}}(M^m) & := & {\displaystyle 
\frac{1}{m-1}\Theta_{f_{g(v)}}(M)+\Pi_{f_{g(v)}}(M)}  \, ,   
\end{array}
\end{equation}
while if $m=1$, 
\begin{equation} \label{eq12b}
\mc{W}_{f_{g(v)}}(M^{1}) := {\displaystyle \Psi_{f_{g(v)}}(M) \equiv 
\Pi_{f_{g(v)}}(M) =: \mc{D}_{f_{g(v)}}(M^{1}) }\, ,
\end{equation}
respectively. In contrast with $\Theta_{f_{g(v)}}$, $\Phi_{f_{g(v)}}$ and
$\Pi_{f_{g(v)}}$, the functionals $\mc{W}_{f_{g(v)}}$ and $\mc{D}_{f_{g(v)}}$ 
are intrinsically defined in the space of metrics $\mc{M}_{[g(v)]}(M)$. 

We exploit the conformal properties of these functionals to define 
special vector field cycles in the tangent bundle, and analyze their 
deformations under scaling of the vector field. We are motivated by the 
observation that Lemma \ref{le2} 
describes the absolute minimum of $\Psi_{f_{g(v)}}(M)$ under conformal 
deformations of the embedding, but $\Psi_{f_{g(v)}}(M)$ 
could admit other extremal 
points that are also critical points of $\mc{W}_{f_{g(v)}}$ under  
conformal deformations of $g_S$ \cite[Theorem 3.10]{gracie}. 

A deformation $v_t=v+w_t$ of $v$ is said to be conformal if
$w_t$ is a time dependent conformal vector field that vanishes
at $t=0$. If 
$\eta_{t}^{w_t}(p)$ is the position at time $t$ along an integral curve of 
$w_t$ that starts at $p$ when $t=0$, then 
$\eta_{t}^{w_t\, *} g = e^{2u_t}g$ for an associated path of functions
$t \rightarrow u_t$ on $M$, with $u_t\mid_{t=0}=0$. In a 
neighborhood of $f_g(v)(M)$, the background Sasaki metric $g_S$ is changing
conformally to $e^{2u_t\circ \pi}g_S$ on $TM$,
and by restriction to $f_{g(v)}(M)$, the metric $g(v_t)$ on the embedded 
submanifold $f(v_t)(M)$ is just $g(v_t)= (e^{2u_t \circ \pi} g)(v)$.
Conversely, a variation of this type defined in a 
neighborhood of $f_g(v)(M)$ in $(TM,g_S)$ arises from a deformation 
of $v$ by a time dependent conformal vector field. We say that $v_t$ is a 
conformal deformation of $v$ of conformal factor $e^{2u_t}$. For notational
convenience, we shall then set
$u_t^v=u_t \circ \pi$.

If $m \neq 2$, we say that $v$ defines a {\it canonical cycle in the homology 
class $[f_{g(v)}(M)]$} if $v$ is a critical point of $\mc{W}_{f_{g(v)}}$ in the 
space of all of its volume preserving conformal deformations, and say that 
$v$ defines a {\it canonical cycle} $f_{g(v)}(M)$ in $(TM,g_S)$ if it is an 
absolute minimizer of $\mc{W}_{f_{g(v)}}$ in the space of all of its conformal 
deformations. Notice that under a conformal deformation $v_t$ of $v$ of 
conformal factor $e^{2u_t}$, if $m>2$ we have that
\begin{equation} \label{eq13}
\begin{array}{rcl}
\mc{W}_{f_{g(v_t)}}(M) \! \! \! \! \! \! &=& \! \! \! \! \! \!{\displaystyle
\int_{f_{g(v)}(M)} \hspace{-0.85cm}  e^{mu^v_t} \left( e^{-2u^v_t} \!
\left( \! \frac{m}{m-1}\! \! \sum_{i,j}K^{TM,g_S}(e^{f_{g(v)}}_i,
e^{f_{g(v)}}_j)\!+ 
\! \| H^{TM,g_S}_{f_{g(v)}}\|^2 \!+m (m\!-\! 2)\| du_t^{v \, \tau}\|^2 
\right) \! \! \right)
d\mu_{f_{g(v)}} }\, , \vspace{1mm} \\
\mc{D}_{f_{g(v_t)}}(M) \! \! \! \! \! \! &= & \! \! \! \! \! \! { \displaystyle
  \int_{f_{g(v)}(M)} \hspace{-0.85cm} e^{mu^v_t}\left( e^{-2u^v_t} \!
\left( \! \frac{1}{m-1}\! \! \sum_{i,j}K^{TM,g_S}(e^{f_{g(v)}}_i,e^{
f_{g(v)}}_j)\!+ 
\! \| \alpha^{TM,g_S}_{f_{g(v)}}\|^2 \! +\! (m\!-\! 2)\| du_t^{v\, \tau}\|^2 
\right) \! \! \right) d\mu_{f_{g(v)}} }\, ,
\end{array}
\end{equation}
respectively, while if $m=1$, all deformations of $v$ are conformal, all 
metrics are flat, and we then have that 
\begin{equation} \label{eq17}
\mc{W}_{f_{g(v_t)}}(M^1) = {\displaystyle
\int_{f_{g(v)}(M)} e^{u^v_t} \left( e^{-2u^v _t}\| H^{TM,g_S}_{f_{g(v)}} \|^2 
\right) d\mu_{f_{g(v)}} }=\mc{D}_{f_{g(v_t)}}(M^1)\, .
\end{equation}
 
If $m=2$ instead, 
$\mc{W}_{f_{g(v_t)}}(M)=\mc{W}_{f_{g(v)}}(M)$, $\mc{D}_{f_{g(v_t)}}(M)=
\mc{D}_{f_{g(v)}}(M)$ and
$S_{f_{g(v_t)}}(M)=\chi(M)$ are all three independent of $t$, so the 
definition 
above requires a refinement. We then consider the functionals 
\begin{equation} \label{eq12c}
\begin{array}{rcl}
\mc{W}^2_{f_{g(v_t)}}(M) & = & {\displaystyle
\int_{f_{g(v)}(M)} \hspace{-0.55cm}  e^{2u^v_t} \left( e^{-2u^v_t} \!
\left( 2\sum_{i,j}K^{TM,g_S}(e^{f_{g(v)}}_i,e^{f_{g(v)}}_j)+ 
\| H^{TM,g_S}_{f_{g(v)}}\|^2 \right) \right)^2 d\mu_{f_{g(v)}} } \, ,
\vspace{1mm}\\ 
\mc{D}^2_{f_{g(v_t)}}(M) & = & {\displaystyle
\int_{f_{g(v)}(M)} \hspace{-0.55cm}  e^{2u^v_t} \left( e^{-2u^v_t} \!
\left( \sum_{i,j}K^{TM,g_S}(e^{f_{g(v)}}_i,e^{f_{g(v)}}_j)+ 
\| \alpha^{TM,g_S}_{f_{g(v)}}\|^2 \right) \right)^2 d\mu_{f_{g(v)}} } \, ,  
\end{array}
\end{equation}
and say that $v$ defines a canonical cycle in the homology class 
$[f_g(v)(M)]$ 
if $v$ is a critical point of $\mc{W}^2_{f_{g(v)}}(M)$ in the space of area 
preserving conformal deformations, while say that $v$ defines a canonical cycle
 $f_{g(v)}(M)$ in $(TM,g_S)$ if it is a minimizer of this functional in the 
space of all conformal deformations.

\begin{theorem} \label{th5}
{\rm a)} A vector field $v$ on $(M^{m=1},g)$ defines a canonical cycle in 
$[f_{g(v)}(M)]$ if, and only if, $v$ is minimal and 
$f_{g(v)}(M)$ is a geodesic in $(TM,g_S)$, and $v$ is a canonical cycle
$f_{g(v)}(M)$ in $(TM,g_S)$ if, and only if, there is no geodesic shorter than
$f_{g(v)}(M)$ in $(TM,g_S)$. If we realize $(M,g)$ as a
geodesic of length $l=\mu_g(M)$ in a flat Hopf torus in $\mb{S}^3$,
then $(TM,g_S)$ is isometric to $\mb{S}^1(l)\times \mb{R}$ with the product
metric, and any of these canonical cycles is given by a vector field embedding 
$M$ into $TM$ as the intersection of $TM$ with a line  
of slope $\frac{2\pi}{l}\mb{Q}$ transversal to the vertical.

\noindent {\rm b)} A vector field $v$ on $(M^{m\geq 2},g)$ defines a canonical 
cycle in $[f_{g(v)}(M)]$ if, and only if, 
$$
\frac{m}{m-1}\sum K^{TM,g_S}( e^{f_{g(v)}}_i,e^{f_{g(v)}}_j)+\| H^{TM,g_S}_{
f_{g(v)}}\|^2 = {\displaystyle \frac{1}{\mu_{g(v)}} \mc{W}_{f_{g(v)}}(M)}\, .
$$
Such a $v$ defines a canonical cycle $f_g(v)(M)$ in $(TM,g_S)$ if in addition 
we have that $H^{TM,g_S}_{f_{g(v)}}=0$.
\end{theorem}

{\it Proof}. a) For dimensional reasons, the curvature of the Sasaki metric 
$g_S$ on $TM$ is zero. We choose a positively oriented frame $\{ T, \nu_H\}$ 
for $TM$ over points of $f_{g(v)}(M)$, so that $T$ is a positive frame for the 
submanifold, and $\nu_H$ is a frame for its normal bundle such that 
$H^{TM,g_S}_{f_{g(v)}}=h_{f_{g(v)}}\nu_H$, where $h_{f_{g(v)}}$ is a scalar 
function. 
By \cite[Theorem 3.3]{gracie}, the critical point equation of (\ref{eq17}) is 
$$
2 \Delta^{g(v)}h_{f_{g(v)}}=2 h^3_{f_{g(v)}} \, . 
$$
By direct computation of the variation of (\ref{eq17}) along deformations of 
the curve, we conclude that a critical point must be such that 
$h_{f_{g(v)}}$ is constant, and since the only
constant solution of the equation above is $h_{f_{g(v)}} \equiv 0$,
it follows that $f_{g(v)}(M)$ is a geodesic in $(TM,g_S)$. 

By \cite[Proposition 1]{pink}, it is always possible to realize $(M,g)$ as
a geodesic of length $l=\mu_g(M)$ in some Hopf tori in $\mb{S}^3$. The 
remaining portion of the proof is now straightforward, cf.  
\cite[Theorem 5]{barr}. 

b) If $m>2$, we set $N=2m/(m-2)$. By (\ref{eq13}) and (\ref{eq12c}),
we have that
$$
\begin{array}{rcl}
{\displaystyle \frac{d}{dt} \mc{W}^2_{f_{g(v_t)}}(M)\mid_{t=0}} & = & 
-{\displaystyle \int_{f_{g(v)}(M)} 2  e^{-2u_t^v}\dot{u}^v_t
\left( 2\sum_{i,j}K^{TM,g_S}(e^{f_{g(v)}}_i,e^{f_{g(v)}}_j)+ 
\| H^{TM,g_S}_{f_{g(v)}}\|^2 \right)^2 d\mu_{f_{g(v)}} } \, , \\
{\displaystyle \frac{d}{dt} \frac{1}{\mu^{\frac{2}{N}}_{g(v_t)}}
\mc{W}_{f_{g(v_t)}}(M)} & = & 
{\displaystyle  \frac{m-2}{\mu^{\frac{2}{N}}_{g(v_t)}} \left(
\int_{f_{g(v)}(M)} \hspace{-0.2cm} e^{mu^v_t}\dot{u^v_t}\left( e^{-2u^v_t}
\left( \frac{m}{m-1} \sum K^{TM,g_S}(e_i^{f_{g(v)}},e_j^{f_{g(v)}}) + 
\| H^{TM,g_S}_{f_{g(v)}} \|^2
\right.  \right. \right. } \vspace{1mm}\\ & & 
-{\displaystyle \left.  \left. \left. 
m(m-2)g(du^{v\, \tau}_t, du^{v\, \tau}_t) + 2m  \Delta^{g(v_t)} u^v_t \right)  -
\frac{1}{\mu_{g(v_t)}}
\mc{W}_{f_{g(v_t)}}(M)\right) \right)  d\mu_{f_{g(v)}} 
 }\, , 
\end{array}
$$
and the vanishing of these expressions for all possible 
$\dot{u}_t$s orthogonal to the constants at $t=0$ implies the existence of
a constant $\lambda_{\mc{W}}$ such that, when $m=2$ we have that
$$
\left( 2\sum_{i,j}K^{TM,g_S}(e^{f_{g(v)}}_i,e^{f_{g(v)}}_j)+ 
\| H^{TM,g_S}_{f_{g(v)}}\|^2 \right)^2 =\lambda_{\mc{W}}^2 \, , 
$$
while when $m>2$, we have that
$$
\frac{m}{m-1}\sum K^{TM,g_S}( e^{f_{g(v)}}_i,e^{f_{g(v)}}_j)+\| H^{TM,g_S}_{
f_{g(v)}}\|^2 = \lambda_{\mc{W}} \, . 
$$
In either case,
$$
\lambda_{\mc{W}}=\frac{1}{\mu_{g(v)}} \mc{W}_{f_{g(v)}}(M)\, .
$$
Thus, $v$ defines a canonical cycle in $[f_{g(v)}(M)]$ if, and only if,
the density of the functional $\mc{W}_{f_{g(v)}}$ equals its average.

Suppose now that $v$ defines a canonical cycle in $(TM,g_S)$. We consider
the scale invariant functionals $\mu_{g(v)}\mc{W}^2_{f_{g(v)}}$ and 
$\mc{W}_{f_{g(v)}}/\mu_{g(v)}^{\frac{2}{N}}$ when $m=2$ and $m>2$, 
respectively. By the partial result above, the vector field $v$ extremizes 
these functionals in the space of all of the gradient deformations of $v$,
and this $v$ is then also a critical point of the volume functional 
(\ref{vgv}), so by (\ref{veqv}), we conclude that 
$\| H^{TM,g_S}_{f_{g(v)}}\|^2=0$.  
\qed

\begin{remark}
If $v$ defines a canonical cycle in the homology class $[f_{g(v)}(M)]$, by
(\ref{eq11}), the metric $g(v)$ will be of constant scalar curvature 
if $v$ extremizes the functional $\mc{D}_{f_{g(v)}}$ also. In fact:
\begin{itemize}
\item If $m\geq 3$, $g(v)$ is a metric of constant scalar curvature
if, and only if, $v$ is also a critical point of $\mc{D}_{f_g(v)}$ in the space 
of volume preserving conformal deformations, and so 
$$
\frac{1}{m-1}\sum_{i,j}K^{TM,g_S}(e^{f_{g(v)}}_i,e^{f_{g(v)}}_j)+ 
\| \alpha^{TM,g_S}_{f_{g(v)}}\|^2 
$$
is constant equal to ${\displaystyle 
\frac{1}{\mu_{g(v)}}\mc{D}_{f_{g(v)}}(M)}$. Among all of these, a $v$ that
defines a canonical cycle $f_{g(v)}(M)$ in $(TM,g_S)$ is a Yamabe metric in its 
conformal class if $f_{g(v)}(M)$ is also a minimizer of
 $\mc{D}_{f_{g(v)}}(M)$ in the said space.
\item If $m=2$, the scalar curvature $s_{g(v)}$ of $g(v)$ is constant if, and 
only if, 
$$
\sum_{i,j}K^{TM,g_S}(e^{f_{g(v)}}_i,e^{f_{g(v)}}_j)+ 
\| \alpha^{TM,g_S}_{f_{g(v)}}\|^2 
$$
is constant and equal to ${\displaystyle \frac{1}{\mu_{g(v)}}
\mc{D}_{f_{g(v)}}(M)}$, and among 
the set of all of these $v$s, a $v$ that defines a canonical cycle 
$f_{g(v)}(M)$ in $(TM,g_S)$ is one for which the value of
$s_{g(v)}^2$ is the smallest.
\end{itemize}
\end{remark}

\begin{corollary} \label{co7}
The zero section $f_{g(0)}(M)\hookrightarrow (TM,g_S)$ is a canonical cycle if,
and only if, the scalar curvature $s_g$ of $g$ is constant.
\end{corollary}

{\it Proof}. At points $(p,0)$ of the zero section, 
$g_S(R_{p,0}^{TM,g_S}(X^{hor},Y^{hor})Z^{hor},W^{hor})=g(R_p^g(X,Y)Z,W)$.
\qed

Corollary \ref{co7} is an expected result. We can always 
find an isometric embedding $f_g: (M^m,g) \rightarrow (\mb{S}^{\tn},\tg)$ into 
a standard sphere of sufficiently large dimension $\tn=\tn(m)$, and cast $TM$ 
as a summand in the Whitney sum decomposition 
$T\mb{S}^{\tn}\mid_{f_g(M)}=
TM\oplus \nu(M)$ of the pullback tangent bundle of this background
sphere over the embedded $M$. As the zero section $f_{g(0)}: 
(M,g) \hookrightarrow (TM,g_S)$, $(M,g)$ has exterior scalar curvature $s_g$
and null mean curvature vector $H^{TM,g_S}_{f_{g(0)}}$, while as the zero 
section of $\nu(M)$, which is identified with a tubular neighborhood of it in
the background sphere, 
it has exterior scalar curvature $m(m-1)$ and mean curvature vector $H_{f_g}$
that when $m\geq 2$, by the existence of isometric embeddings $f_{g_t}$ of 
volume preserving conformal metric deformations of $g$, can be proved to be 
of constant norm \cite[Theorem 6]{scs}, and so
$s_g$ is constant if, and only if, the second fundamental form $\alpha_{f_g}$ 
of $f_g$ has constant norm. This is the 
essence of the proof of Theorem \ref{th5} we provided for manifolds of 
dimension 1, where our statements are all transparent since $g$ and 
$g_S$ are then necessarily flat, and in which case the geodesic zero section of
$(TM,g_S)$ is also a geodesic zero section of $\nu(M)$ if, and only if,
the length $\mu_g(M)=2\pi$.   

Suppose that $v$ is a nontrivial vector field on 
$(M^m,g)$, whose singular set $Z(v)$, if not empty, is a compact submanifold 
of $M$. When defined,  
we denote by $N_{\varepsilon}(Z(v))$ the $\varepsilon$ tubular neighborhood of 
$Z(v)$. Then if $t\in \mb{R}_{>0}$, the $m+1$ dimensional sheet in $(TM,g_S)$ 
given by  
$$
M^{\varepsilon}_{v,t}=\{(p,st(v/\| v\|_g)(p)): \; p\in M\setminus 
N_{\varepsilon} (Z(v))\, , \; 0\leq s \leq 1\} \, , 
$$  
is a rectifiable manifold with rectifiable boundary with edges or corners. If 
$\iota $ is the inclusion map, we define the $m+1$ sheat current 
$[[ M_{v,t}]]$ as the weak* limit 
$$
C^{\infty}(TM;\Lambda^{m+1} TM) \ni \omega \rightarrow 
[[ M_{v,t}]](\omega) := 
\lim_{\varepsilon \searrow 0} \int_{M^{\varepsilon}_{v,t}}\iota^* \omega 
\, .
$$
The homological boundary of this current is the difference of a suitable 
current on $f_{g(tv/\|v\|_g)}(M\setminus Z(v))$, and the current of 
integration over the zero section $f_{g(0)}(M)$, and so when acting on a 
closed 
form $\alpha$, this difference is zero, and the summands encode the   
topological information given by evaluating the cohomology class density 
$[f^*_{g(0)}(\alpha)]$ over $M$, and how to compute it directly from 
$[\alpha]$, $v$ and $Z(v)$, by a current supported in a neighborhood of 
$f_{tv/\| v\|_g}(M\setminus Z(v))$ in $TM$.   

\section{Some illustrative examples and related results}

\subsection{The standard sphere $(\mb{S}^{2n+1},g) \hookrightarrow 
(\mb{R}^{2n+2}, \| \; \|^2)$} 
We consider firstly the odd dimensional unit sphere
$$
\mb{S}^{2n+1}=\{ z=(z_0, \ldots , z_{n}) \in \mb{C}^{n+1}\equiv \mb{R}^{2n+2}: 
\; |z|=1 \} 
$$
as a hypersurface in the complex vector space $\mb{C}^{n+1}$. 
If $z_k=x_{2k}+ix_{2k+1}$, the vector fields $\{H_k = x_{2k}\partial_{
x_{2k+1}}- x_{2k+1}\partial_{x_{2k}}\}_{k=0}^n$ form a basis of the Lie 
algebra $\mf{t}_{n+1}$ of a maximal torus $T_{n+1}$ in the automorphism group 
$\mb{U}(n+1)$ of its standard almost contact structure of 
contact form $\sigma= \sum_{k=0}^n (x_{2k}dx_{2k+1}-x_{2k+1}
dx_{2k})$, Reeb vector field $v_H= \sum_{k=0}^n H_k$, and $(1,1)$-tensor 
$\phi$ defined by the restriction of the complex structure $J$ of 
$\mb{C}^{n+1}$ on $\mc{D}={\rm ker}\, \sigma$, and $\phi (v_H)=0$, 
respectively. The ensuing compatible metric 
$$
g=d\sigma \circ (\BOne \otimes \phi) + \sigma \otimes \sigma 
$$
is the standard (Einstein) metric on $\mb{S}^{2n+1}$, of scalar curvature
$s_g=2n(2n+1)$. Thus, the almost contact structure $(v_H,\sigma,\phi)$ is 
normal, and the compatible metric structure $g$ is transversally extremal. 
The mapping $\mb{R}^{n+1}\ni w=(w_0, \ldots, w_n) \mapsto v_{H_w}=
\sum_k w_k H_k \in \mf{t}_{n+1}$ identifies the Lie algebra $\mf{t}_{n+1}$ of 
the maximal torus $T_{n+1} \subset \mb{U}(n+1)$ with $\mb{R}^{n+1}$. 

We consider now the set $\mc{S}(\mc{D},J)$ of all compatible contact metric 
structures associated with the strictly pseudoconvex CR structure $(\mc{D},J)$
that $\sigma$ and $J$ define. A vector field $X$
that preserves $(\mc{D},J)$ is said to be positive if 
$\sigma(X)>0$, and any such is conjugate to a $v_{H_w} \in \mf{t}_{n+1}$ with 
$w_i >0$, thus identifying the set of all of them with $\mb{R}_{>0}^{n+1}$. The
foliation that each positive $v_{H_w}$ defines is taut, 
and all of these fields yield almost contact metric structures that can be 
represented by 
a compatible transversally extremal contact metric \cite[Example 7.7]{bgsi}, 
with those associated to the ray defined by $v_H$ itself being the only ones 
of constant scalar curvature.

Indeed, if $w\in \mb{R}_{>0}^{n+1}$, then  
$(v_{H_w},\sigma_w:=\sigma/\sigma(v_{H_w}),\phi_w,g_w)$,
where $\phi_w$ is defined by $J\mid_{\mc{D}}$ on $\mc{D}$
and $\phi_w(v_w)=0$, respectively, and $g_w$ is given by 
\begin{equation} \label{cme}
g_w=d\sigma_w \circ (\BOne \otimes \phi_w) + \sigma_w \otimes \sigma_w \, , 
\end{equation}
is a contact metric structure, and the leaves of the foliation of $v_{H_w}$ 
are $g_w$ geodesics. If 
$w \in \mb{Z}_{>0}^{n+1}$, $(v_{H_w},\sigma_{H_w},\phi_w,g_w)$ 
is quasi-regular, and its transversal is the orbifold weighted 
projective space $(\mb{P}^n_w(\mb{C}),J_w)$. The space of metrics 
$\mc{M}(v_{H_w},J_w)$ associated with the polarized manifold 
$(\mb{S}^{2n+1},v_{H_w},J_w)$ has a representative $\tilde{g}_w$ of
transverse K\"ahler metric $\tilde{g}_w^T$ of scalar curvature  
$$
s_{\tilde{g}_w^T}=4(n+1)\frac{\sum_{j=0}^n w_j(2(\sum_{k=0}^n w_k) -(n+2)w_j)
|z_j|^2 }{\sum_{j=0}^n w_j |z_j|^2 }\, ,
$$
and mean transverse scalar curvature $4n\sum _{j=0}^n w_j$. We have that
$$
\mu_{\tilde{g}_w}(\mb{S}^{2n+1})=2\frac{\pi^{n+1}}{n!}\frac{1}{\prod_{j=0}^n 
w_j} \, ,
$$
and since $s_{\tilde{g}_w}=s_{\tilde{g}_w^T}-2n$, the projection 
$s^0_{\tilde{g}_w}$ 
of $s_{\tilde{g}_w}$ 
onto the constants is given by
$$
s^0_{\tilde{g}_w}=2n(2 \sum_{j=0}^n w_j -1)\, , 
$$
so
$$
s_{\tilde{g}_w}-s^0_{\tilde{g}_w}=4(n+2)\frac{\sum_{j=0}^n w_j((\sum_{k=0}^n 
w_k) -(n+1)w_j) |z_j|^2 }{\sum_{j=0}^n w_j |z_j|^2 }\, .
$$
If $X=\partial_{J_w}^{\#}f$ where 
$f=\sum_{i=0}^n b_i |z_i|^2$ is a basic real holomorphy potential, the 
Sasaki-Futaki character of the polarization $(v_{H_w},J_w)$ acting on $X$ is  
\begin{equation} \label{SF}
\begin{array}{rcl}
\mf{F}_{(v_{H_w},J_w)}(X) & = &{\displaystyle -8(n+2)\pi ^{n+1} 
\int_{\mb{R}^{n}_{+}}
\frac{(b_0+\sum_{j=1}^n b_jx_j)(w_0A_0+\sum_{j=1}^n w_jA_jx_j )}
{(w_0+\sum_{j=1}^n w_j x_j)^{n+3} }dx_1\ldots dx_n } \vspace{1mm} \\
& = & {\displaystyle -16\frac{\pi ^{n+1}}{(n+1)!}\left(
\sum_{i=0}^n \frac{b_i}{w_i}A_i +\frac{1}{2}\sum_{i\neq j}\frac{b_i}{w_i}A_j
\right) \frac{1}{\prod_{j=0}^n w_j}}\, ,
\end{array}
\end{equation}
where, for convenience, we have set $A_j=\sum_{k=0}^n w_k -(n+1)w_j$. If
$w\in \mb{Q}^{n+1}_{>0}$ instead, we can use a contact homothety 
dilating $w$ into a weight in $\mb{Z}^{n+1}_{>0}$, and show the  
existence of a $\tilde{g}_w$ with transversally extremal metric $g^T_w$ from 
the preliminary result above for integer weights. We can then apply the 
openess theorem \cite[Theorem 7.6]{bgsi} to arrive at the same conclusion for 
arbitrary weights $w\in \mb{R}^{n+1}_{>0}$.
Notice that $\mf{F}_{(v_{H_w},J_w)}\equiv 0$ if, and only if,
$w$ is of the form $w=l(1, \ldots, 1)$, and the scalar curvature of 
$\tilde{g}_w^T$ is the constant $s_{\tilde{g}_w^T}=  
4n(n+1)l$, and that if $w=(1,\ldots, 1)$, then $v_{H_w}=v_H$, 
$\tilde{g}_w =g$ and $\tilde{g}_w^T$ is the standard Fubini-Study metric 
$g_{FS}$ on $\mb{P}^n(\mb{C})$ lifted as a metric on the horizontal space
at any point of the fiber.   

The field $v_H \in \mf{t}_{n+1}$ is the only positive element that lies 
on the unit sphere tangent bundle of the standard sphere. It is tangential to 
the geodesic circle fibers of the Hopf Riemannian submersion 
$\mb{S}^1 \hookrightarrow (\mb{S}^{2n+1},g) \rightarrow (\mb{P}^n(\mb{C}),
g_{FS})$, and arises as the infinitesimal generator of the $\mb{S}^1$ 
action on $\mb{S}^{2n+1}$ leading to it, thus named 
the Hopf vector field of $(\mb{S}^{2n+1},g)\hookrightarrow 
(\mb{C}^{n+1},\sum dz_i \otimes d\overline{z}_i)$. It may be cast explicitly 
in 
terms of the CR structure $(\mc{D},J)$. For the normality of the almost 
contact structure  $(v_H,\sigma,\phi)$ may be phrased in terms of the cone 
$C(\mb{S}^{2n+1})=\mb{R}_{>0}\times \mb{S}^{2n+1}$ provided with the metric
$$
g_C=dr^2 + r^2 g\, .
$$
If $r\partial_r$ is the radial vector field, then 
$\mc{L}_{r\partial_r}g_C=2g_C$, and the cone carries the almost complex 
structure $I_J$ defined by
$$
I_J(Y)=\phi(Y)-\sigma(Y)r\partial_r \quad {\rm and}\quad 
I_J(r\partial_r )=v_H \, , 
$$
respectively. The said normality is equivalent to the integrability of $I_J$.   
Notice that $\mu_{g(v_H)}^{T\mb{S}^{2n+1}}(v_H)==2^n \mu_g(\mb{S}^{2n+1})$, 
and that
except when $n=0$, $v_H$ is not parallel and $[g(v_H)] \neq [g(0)=g]$. 

A geodesic embedding $\iota_{2n>0}: \mb{S}^{2n} \hookrightarrow 
\mb{S}^{2n+1}$ produces a Whitney sum decomposition $\iota_{2n}^* 
T\mb{S}^{2n+1}=T\mb{S}^{2n} \oplus \nu(\mb{S}^{2n})$, with a well-defined
projection $\pi^{2n}$ onto the first factor. We can then define 
the vector field $v_{H,\iota}^{2n}$ of $\mb{S}^{2n}$ by 
$v_{H,\iota}^{2n}=\pi^{2n}_* (\iota_{2n}^* v_H)$, and since the complex
structure $J$ fixes the CR structure on the equatorial spheres 
$(\mb{S}^{2n-1},g) \hookrightarrow (\mb{C}^{n+1},dz\otimes d \overline{z})$    
in $(\mb{S}^{2n+1},g)$ defined by the vanishing of some complex coordinate
$z_k$, up to a rotation on this $z_k$ plane, we define the Hopf vector field 
$v_{H}^{2n}$ by this expression when $\iota$ is the geodesic embedding
of $\mb{S}^{2n}$ given by ${\rm Im}(z_k)=0$.
The zero set $Z(v_H^{2n})$ 
consists of antipodal points $\{N,S\}$ of $\mb{S}^{2n}$ defined by
${\rm Re}(z_k) = \pm 1$.  
We set $[[ f_{g(t v_H^{2n}/\| v_H^{2n}\|_g)}(\mb{S}^{2n}\setminus \{ N,S\})]]$ 
to be the weak* limit as $\varepsilon \searrow 0$ of the $2n$ 
current of integration over $f_{g(tv_H^{2n}/\|v_H^{2n}\|_g)}(\mb{S}^{2n}
\setminus N_{\varepsilon}(Z(v_H^{2n}))$ in $(TM,g_S)$, 
while $[[ f_{g(0)}(\mb{S}^{2n})]]$ is the $2n$ current of integration over the
zero section. 

\begin{example} \label{hfs}
Let $v_H$ and $v_H^{2n}$ be the Hopf vector field on $(\mb{S}^{2n+1},g)
\hookrightarrow (\mb{C}^{n+1},\sum dz_i \otimes d\overline{z}_i)$ and
its projection onto $\mb{S}^{2n}$, respectively. Then 
\begin{enumerate}[label={\rm (\alph*)}]
\item {\rm (\cite{glzi,olel,olel2})} $v_H$ is a minimal unit vector field on 
the standard sphere $(\mb{S}^{2n+1},g)$.
\item For any $t\in \mb{R}_{>0}$, $f_g(tv_H)(\mb{S}^{2n+1\geq 3})$ is a 
canonical 
cycle in the homology class $[f_g(tv_H)(\mb{S}^{2n+1})]$, $g(tv_H)$ has
positive sectional curvature, and as $t\searrow 0$, the cycle converges 
uniformly to the 
zero section canonical cycle of $(T\mb{S}^{2n+1},g_S)$, the absolute minimizer 
of $\Psi_{f_g(v)}$.
\item The vector field $v^{2n}_H$ is a minimal unit singular vector field on 
the standard sphere $(\mb{S}^{2n\geq 2},g)$, with zeroes at antipodal points 
$\{N,S\}$. 
\item For any $t \in \mb{R}_{>0}$,  the sectional curvatures of
$f_g(tv^{2n}_H/\|v^{2n}_H\|_g)(\mb{S}^{2n} \setminus \{ N,S\})$ are bounded
above, are strictly positive everywhere if $t<\sqrt{4/3}$, and as 
$t\searrow 0$, they all converge pointwise but not uniformly to $1$.  
As $2n$ currents, the homological boundary of the sheat current 
$[[\mb{S}^{2n}_{v_{H}^{2n},t}]]$ is given by 
$$
\partial [[\mb{S}^{2n}_{v_{H}^{2n},t}]] = 
[[ f_{g(t v_H^{2m}/\| v_H^{2m}\|_g)}(\mb{S}^{2n}\setminus \{ N,S\})]] -
[[ f_{g(0)}(\mb{S}^{2n})]] \, , 
$$
and so if $\alpha \in C^{\infty}(T\mb{S}^{2n}; \Lambda^{2n}T\mb{S}^{2n})$ is
closed, 
$$
[[ f_{g(t v_H^{2m}/\| v_H^{2m}\|_g)}(\mb{S}^{2n}\setminus \{ N,S\})]](\alpha)=
\int_{\mb{S}^{2n}}f_{g(0)}^{*} \alpha \, . 
$$
\end{enumerate}
\end{example}

{\it Proof}. (a) The Hopf vector field $v=v_H=J(\partial_r)$ on 
$(\mb{S}^{2n+1},g) \hookrightarrow (\mb{C}^{n+1}=\mb{R}^{2n+2},\| \, \cdot \, 
\|^2)$ is the Reeb vector field of the strictly 
pseudoconvex CR structure $(\mc{D},J)$, and its flow is by isometries. We 
choose a local orthonormal tangent frame $\{ e_i \}_{i=0}^{2n}$ about 
an arbitrary point $p$, with $e_0=v$ satisfying the globally defined relation 
$\nabla^g_{e_0}v=0$, while the remaining elements of the frame are all 
horizontal, constitute an orthonormal set whose span pointwise is left 
invariant under the action of $J$, and for each one of them, 
$\nabla^g_{e_i} v= Je_i$, and  $C_v^t C_v (e_i)=(\nabla^g v)^t J{e_i}=e_i$, 
respectively. Hence, 
$\lambda_0=0$, $c^{-2}_0=1$, and for $i>0$,  
$\lambda_i =1$ and $c^{-2}_i=2$, correspondingly. The set
$\{ e_0^{ver}, \{ -((\nabla v)^t(e_i))^{hor}+e_i^{ver}\}_{
i=1}^{2n}\}$ is a basis for the normal bundle of 
$f_{g(v)}({\mb{S}^{2n+1}})$ in $T(\mb{S}^{2n+1})$ at $p$, 
while the set $\{-((\nabla v)^t(e_i))^{hor}+e_i^{ver}\}_{
i=1}^{2n}\}$ is an basis of the normal bundle to the submanifold
in $S^1(T\mb{S}^{2n+1})$.

If $i=0$, $\nabla^g_{e_i}e_i$ and $R^g(v,\nabla^g_{e_i}v)e_i$ equal the 
null vector, while if $i>0$, $(\nabla^g_{e_i}e_i)^{ver} = \< e_i, Je_i\> v = 
0$, $\nabla^g_{e_i}e_i\mid_p=0$ and 
$R^g(v,\nabla^g_{e_i}v)e_i= \< \nabla^g_{e_i}v,e_i\> v -\< v,e_i\> 
\nabla^g_{e_i}v= \< Je_i,e_i\>v= 0$ also. 
On the other hand, 
$\nabla^g_{e_0}\nabla^g_{e_0}v=0$, while
if $i>0$, $\nabla^g_{e_i}\nabla^g_{e_i}v=\< e_i, J^2 e_i\> v=-v$.
Thus,
$$
\sum_{i\geq 0} c_i^2\left(\nabla^g_{e_i}\nabla^g_{e_i} -
\nabla^g_{\nabla^g_{e_i} e_i} -\nabla^{g}_{R^g(v,\nabla^g_{e_i} v)e_i}\right)v
=-n v \, , 
$$
and identity (\ref{eq10}) holds with $\lambda =-n$.

(b) We let $\{ e_i\}_{i=0}^{2n}$ be the orthonormal frame in (a) above, and
consider the orthonormal tangent frame 
$\{ e_i^{f_{g(tv_H)}}=c_i(e_i^{hor}+(\nabla^g_{e_i}tv_H)^{ver})\}_{i=0}^{2n}$,
$c_i^{-2}= 1+\| \nabla^g_{e_i}tv_H\|^2$, of the embedded 
manifold $f_{g(tv_H)}(\mb{S}^{2n+1})$   
defined by the scaled Hopf vector field $tv_H$. At any point $(p,tv(p))$ in 
the embedded manifold, the set 
$\{ N_k =c_k(-((\nabla^{g}tv_H)^t(e_k))^{hor}+ e_k^{ver}\}_k$
is now an orthonormal frame of the normal bundle. 
By (\ref{mcv}), 
the mean curvature vector $H^{T\mb{S}^{2n+1},g_S}_{f_{g(tv_H)}}$ is 
orthogonal to
all the $N_k$s for $k>0$, and along $N_0$, it equals $-n(tv_H)$. We thus have
that $\| H^{T\mb{S}^{2n+1},g_S}_{f_{g(tv_H)}}\|^2=t^2n^2$. 

Since $R^g$ is covariantly constant, 
the sectional curvatures of the various sections (spanned by the elements of 
the frame $\{e^{f_{g(tv_H)}}_i\}$) can be found to be given by 
$$
\begin{array}{rcl}
K^{T\mb{S}^{2n+1},g_S}(e^{f_{g(tv_H)}}_0,e^{f_{g(tv_H)}}_{j>0})& = & 
{\displaystyle \frac{1}{1+t^2}\left(1-\frac{3}{4}t^2+ \frac{1}{4}t^4
\right)}\, , \vspace{1mm}\\
K^{T\mb{S}^{2n+1},g_S}(e^{f_{g(tv_H)}}_{i>0},e^{f_{g(tv_H)}}_{i\neq j>0}) & 
= & {\displaystyle \frac{1}{(1+t^2)^2}\left(1+ 2t^2 \<Je_i,e_j\>^2\right)}\, ,
\end{array}
$$
respectively. The exterior scalar curvature on the submanifold 
$f_{g(tv_H)}(\mb{S}^{2n+1})$ is  
$$
\sum_{i,j}K^{T\mb{S}^{2n+1},g_S}(e^{f_{g(tv_H)}}_i,e^{f_{g(tv_H)}}_{j})= 
\frac{1}{(1+t^2)^2}n(4n+2+5t^2-2t^4+t^6) \, .
$$
By Theorem \ref{th5}, $tv_H$ defines a canonical cycle in the homology
class $[f_{g(tv_H)}(\mb{S}^{2n+1})]$, and by the expressions given above, all 
sectional
curvatures of this cycle converge uniformly to $1$ as $t \searrow 0$, so
the canonical cycles $f_{g(tv_H)}$ in the homology classes $[f_g(tv_N)]$ 
converges uniformly to the canonical cycle given by the zero section
totally geodesic submanifold $(\mb{S}^{2n+1},g)$ of $(T\mb{S}^{2n+1},g_S)$. 

(c) We choose the north and south pole of $\mb{S}^{2n+1}$ to be the points
$(0, \ldots, 0, \pm i)$, and use the natural geodesic embedding
$\mb{S}^{2n} \hookrightarrow \mb{S}^{2n+1}$ defined by 
$$
(z',x_n):=(z_0, \ldots, z_{n-1},x_n)\mapsto (z_0, \ldots, z_{n-1},x_{2n}\pm i 
0) \, .
$$ 
Then $v_H^{2n}= \sum_{k=0}^{n-1} H_k$, and 
$\| v_H^{2n}\|^2= \sum_{k=0}^{n-1} |z_k|^2 = 1-x_{2n}^2$, so $v^{2n}_H$  has 
zeroes at $N=(0, \ldots, 0, 1)$ and $S=(0, \ldots, 0,-1)$, the
north and south pole of $\mb{S}^{2n}$. We set $v=\frac{1}{\sqrt{1-x_{2n}^2}}
v_H^{2n}$ and $e_r= \frac{x_{2n}}{\sqrt{1-x_{2n}^2}}\sum_{j=0}^{n-1}(x_{2j} 
\partial_{x_{2j}}+x_{2j+1}\partial_{x_{2j+1}})-\sqrt{1-x_{2n}^2}
\partial_{x_{2n}}$,  
respectively. 

We consider an orthonormal frame $\{ e_0=v, e_1, \ldots, e_{2n-2},e_{2n-1}=
e_r\}$, where for each $x_{2n} \in (-1,1)$, $e_1, \ldots, e_{2n-2}$ is an 
orthonormal frame of horizontal vectors of the scaled Hopf fibration 
over the $2n-1$ sphere of radius $\sqrt{1-x_{2n}^2}$. Since 
$\nabla^g_{e_0}e_0=-\frac{x_{2n}}{\sqrt{1-x_{2n}^2}}e_r$, $\nabla^g_{e_r}e_0=0$,
 $\nabla^g_{e_0}e_r=\frac{x_{2n}}{\sqrt{1-x_{2n}^2}}e_0$, 
$\nabla^g_{e_r}e_r=0$, and as observed in (a), 
$\nabla^g_{e_i}e_0=\frac{1}{\sqrt{1-x_{2n}^2}}Je_i$,
$\nabla^g_{e_i}\nabla^g_{e_i}e_0=-\frac{1}{1-x_{2n}^2}v$,
 and $(\nabla^g_{e_i}e_i)^{ver}=0$, $\nabla^g_{e_i}e_i\mid_p=0$
for $i=1, \ldots, 2(n-1)$, we have that  
$c_0^{-2}=\frac{1}{1-x_{2n}^2}$, $c^{-2}_{2n-1}=1$, 
$c_i^{-2}=1+\frac{1}{1-x_{2n}^2}$ for $i=1, \ldots, 2(n-1)$, and 
$\nabla^g_{e_0} e_0 +R^g(v,\nabla^g_{e_0} v)e_0 =0$. 
It follows that
$$
\sum_{i\geq 0} c_i^2\left(\nabla^g_{e_i}\nabla^g_{e_i} -
\nabla^g_{\nabla^g_{e_i} e_i} -\nabla^{g}_{R^g(v,\nabla^g_{e_i} v)e_i}\right)v 
=\sum_{i\geq 0} c_i^2 \nabla^g_{e_i}\nabla^g_{e_i}v 
= -\left(x_{2n}^2+\frac{2(n-1)}{2-x_{2n}^2} \right) v\, ,    
$$
so identity (\ref{eq10}) holds with 
$\lambda = -\left(x_{2n}^2+\frac{2(n-1)}{2-x_{2n}^2}\right)$.

(d) By (c) above, the squared norm of the mean curvature vector of the 
(open) $m$ manifold 
$f_{g(tv^{2n}_H/\|v^{2n}_H\|_g)}(\mb{S}^{2n}\setminus \{N,S\})$ is 
$\| H^{T(\mb{S}^{2n}\setminus \{N,S\}),g_S}_{f_{g(tv)}}\|^2=
t^2\left( \frac{x_{2n}^2}{1+(t^2-1)x_{2n}^2}+\frac{2(n-1)}{2-x_{2n}^2} \right)^2$, and its various
sectional curvatures are given by 
$$
\begin{array}{rcl}
K^{T\mb{S}^{2n}\setminus \{N,S\},g_S}(e^{f_{g(tv)}}_0,e^{f_{g(tv)}}_{2n-1})& = 
& {\displaystyle \frac{1-x_{2n}^2}{1+(t^2-1)x_{2n}^2}\left(1-\frac{3}{4}t^2+
\frac{x_{2n}^2}{1-x_{2n}^2} \frac{1}{4}t^4\right)}\, ,\vspace{1mm}\\
K^{T\mb{S}^{2n}\setminus \{N,S\},g_S}(e^{f_{g(tv)}}_0,e^{f_{g(tv)}}_{j\neq 
0,2n-1})
& = & {\displaystyle \frac{1-x_{2n}^2}{(1+(t^2-1)x_{2n}^2)}\frac{1-x_{2n}^2}{(1
-x_{2n}^2+t^2)} \left(1-\frac{3}{4}\frac{t^2}{1-x_{2n}^2}+ \frac{1}{4}
\frac{t^4}{(1-x_{2n}^2)^2} \right)}\, ,\vspace{1mm}\\
K^{T\mb{S}^{2n}\setminus \{N,S\},g_S}(e^{f_{g(tv)}}_{i\neq 0,2n-1},
e^{f_{g(tv)}}_{j\neq 0,i, 2n-1}) & = &
{\displaystyle \left( \frac{1-x_{2n}^2}{1-x_{2n}^2+t^2} \right)^2
\left(1+ 2 \frac{t^2}{1-x_{2n}^2} \<Je_i,e_j\>^2\right)}\, , \vspace{1mm} \\
K^{T\mb{S}^{2n}\setminus\{N,S\},g_S}(e^{f_{g(tv)}}_{i\neq 0,2n-1},
e^{f_g(tv)}_{2n-1}) & = & {\displaystyle \frac{1-x_n^2}{1-x_{2n}^2+t^2} }\, ,
\end{array}
$$
respectively. The cycle defined by $tv_{H}^{2n}$ is not canonical, and
the statements about is sectional curvatures are straightforward  
conclusions from the sectional curvature expressions above.

We consider the $2n+1$ sheet $\mb{S}^{2n \hspace{1mm} \varepsilon}_{v_{H}^{2n},
t}$, and apply Fubini's theorem in the vertical direction to write 
$[[{\mb{S}^{2n \hspace{1mm} \varepsilon}}_{v_{H}^{2n},t}]](d\alpha)$  
as the integration of $\alpha$ over the $2n$ chain sum of the  
three type of components of $\partial \mb{S}^{2n \hspace{1mm} 
\varepsilon}_{v_{H}^{2n},t}$. The measure of the cylindrical sides  of 
this chain is of the order $\varepsilon^{2n-1}$, so the current of integration
of $\alpha$ over it goes to zero as $\varepsilon \searrow 0$. The
desired result follows.   
\qed    
 
As we indicated earlier, by modifying compatibly the Riemannian structure 
on the sphere, all the vector fields $\pm v_{H_w}\in \pm \mb{R}^{n+1}_{>0}$ 
can be made Killing fields of the modified metric. In fact, 
they are minimal unit. In complete generality, we have the following. 

\begin{example}
On $(\mb{S}^{2n+1},g^{v_{H_w}})$, where $g^{v_{H_w}}$  
is any metric compatible with an almost contact metric structure of contact 
bundle $(\mc{D},J)$ and Reeb vector field $v_{H_w}$, the vector fields
$\pm v_{H_w}$ are minimal unit, and the zero section of 
$(T\mb{S}^{2n+1},g_{\, S}^{v_{H_w}})$ is a canonical cycle if, and only if, 
$w=l(1, \ldots, 1)$ and $g^{v_{H_w}}$ is transversally extremal.
\end{example}

{\it Proof}. We consider the metric $g_w$ in (\ref{cme}) compatible with the 
almost contact structure $(v_{H_w},\sigma_w, \phi_w)$ thereby defined. The  
set of all quadruples $(v_{H_w}, \tilde{\sigma}, \tilde{\phi},\tilde{g})$ of 
almost contact metric structures of contact bundle $(\mc{D},J)$ and Reeb field
$v_{H_w}$ correspond to the same basic cohomology class 
$[d\tilde{\sigma}]_B=[d\sigma_w]_B$ and transverse holomorphic structure $J_w$
on $\mb{P}_{w}(\mb{C})$, so in its entirety, it is an affine space modeled on
the product $C^{\infty}_B(\mb{S}^{2n+1}) \times C^{\infty}_B(\mb{S}^{2n+1})
\times H^1(\mb{S}^{2n+1},\mb{Z})=C^{\infty}_B(\mb{S}^{2n+1}) \times 
C^{\infty}_B(\mb{S}^{2n+1})$, where
$C^{\infty}_B(\mb{S}^{2n+1})$ is the space of basic functions, 
 and as we indicated earlier, we can always deform 
$g_w$ to a $\tilde{g}$ in this space so that the transversal K\"ahler metric
of $\tilde{g}$ is extremal, with transversal K\"ahler form 
$d\tilde{\sigma}=d\sigma_w + 
i \partial_{J_w} \overline {\partial }_{J_w} f$ for some basic function $f$.

Given the metric $\tilde{g}$ of any almost contact metric structure 
$(v_{H_w}, \tilde{\sigma}, \tilde{\phi},\tilde{g})$ as above, 
the Riemannian submersion of the fibration onto 
the weighted projective space has totally geodesic fibers, and its O'Neil 
tensor $A^w$ is given by $A_XY=-d\tilde{\sigma}(X,Y) v_{H_w}$ and 
$A^w_X v_{H_w}= J_w X$, where $X,Y$ are horizontal. About any point on the
manifold, we can choose an orthonormal frame $\{ e_j\}_{j=0}^{2n}$ with 
$e_0=v_{H_w}$ and $e_1, \ldots, e_{2n}$ horizontal vectors. Then 
$\nabla^{\tilde{g}}_{e_0}v_{H_w}=0$, $\nabla^{\tilde{g}}_{e_i}v_{H_w}=J_w e_i$, 
and $(\nabla^{\tilde{g}}_{e_j}e_j)^{ver}=0$ for all $j\neq 0$, so  
 $c_j^{-2}=1+ \|\nabla^{\tilde{g}}_{e_j}v_{H_w}\|^2$ equals $1$ or $2$ if $j=0$
or $j=1, \ldots, 2n$, respectively, and by the identity 
$R^{\tilde{g}}(X,v_{H_{w}})Y=\tilde{\sigma}(Y)X-\tilde{g}(X,Y)v_{H_w}$, 
we obtain that $R^{\tilde{g}}(v_{H_{w}},\nabla^{\tilde{g}}_{e_j} v_{H_w})e_j=0$
for all $j$. Thus, (\ref{eq10}) holds with $\lambda = -n$.  

The vanishing of the character (\ref{SF}) of the polarization 
$(v_{H_w},J_w)$ is a necessary condition for the transversal metric 
of $\tilde{g}$ to be extremal of constant scalar curvature. The last 
statement now is a direct application of Corollary \ref{co7}.
\qed     

On the other hand, we can consider $(\mb{S}^{2n+1},g)\hookrightarrow 
(\mb{R}^{2n+2},\| \; \|^2)$ fixing its metric, and take an arbitray compatible
complex structure $I$ in the Euclidean space. This $I$ then induces 
a strictly pseudoconvex CR structure on $\mb{S}^{2n+1}$ of contact
type and Reeb vector field $I(\pm r\partial_r \mid_{\mb{S}^{2n+1}})$. A vector 
field of this type on $(\mb{S}^{2n+1},g)$ is called the Hopf vector field of 
the compatible complex structure $I$ in Euclidean space used to define it, and 
they are all minimal unit of volume $2^2\mu_g(\mb{S}^{2n+1})$ 
\cite[Definition 3.1, Proposition 3.2]{ol}. 
The strictly pseudoconvex CR structure on $\mb{S}^{2n+1}$ depends upon the
complex deformation class of $I$ in $(\mb{R}^{2n+2},\| \; \|^2)$, and remains
fixed when the complex deformation is trivial, (that is to say, 
a complex deformation that infinitesimally is of the form $\mc{L}_X I$ for 
a vector field $X$ whose flow generates elements of $\mb{S}\mb{O}(2n+2)$,) and
since $H^1(\mb{C}^{n+1}, \Theta)$ vanishes, this CR structure is rigid.

\begin{theorem}
If $v$ is a minimal unit vector field on $(\mb{S}^{2n+1},g)\hookrightarrow
(\mb{R}^{2n+2}, \| \;  \|^2)$, then $v$ is, up to a conjugating diffeomorphism,
the Hopf vector field of a compatible complex structure $I$ in $\mb{R}^{2n+2}$,
 $\tilde{v}=g(v, \, \; \, )$ is a contact form, and   
$({\rm ker} \, \tilde{v}, I\mid_{{\rm ker} \tilde{v}})$ is a contact 
strictly pseudoconvex {\rm CR} structure on $\mb{S}^{2n+1}$ with Levi form 
$$
-\tilde{v}([X,IY])=\omega_v(X,IY)=g(X,Y) 
$$
positive or negative depending upon the orientation.
\end{theorem}

{\it Proof}. By Theorem \ref{th4}, $v$ is a Killing field, and $g$ is
bundle like relative to the the flow. Since the geodesics on the sphere are
great circles, the geodesible flow of $v$ helps define a free and proper 
$\mb{S}^1$ action on $\mb{S}^{2n+1}$. We let $\mb{P}_v$ be the manifold of 
orbits, an use the projection map to obtain a Riemannian submersion 
$(\mb{S}^{2n+1},g)\rightarrow (\mb{P_v,\pi_*} g)$, with totally geodesic
$\mb{S}^1$ fibers.    

Along geodesic segments, $v$ restricts to a Jacobi field. Hence, if
$X$ is a normal vector orthogonal to $v$ at $p$, and $\gamma_X$ is a 
geodesic through $p$ with initial velocity $X$, we have that
$$
-v=R^g(\dot{\gamma}_X , v)\dot{\gamma}_X\mid_{t=0}=\nabla^g_{\dot{\gamma}_X}
\nabla^g_{\dot{\gamma}_X }v \mid_{t=0} = \nabla^g_{X} \nabla^g_X v\, , 
$$
so the map
$$
\mc{H} \ni X \mapsto I_vX:=\nabla^g_X v \in \mc{H} 
$$
defines a norm preserving tensorial automorphism of the hyperplane 
distribution $\mc{H}$ of horizontal vectors orthogonal to $v$, and since 
$\nabla^g v=\omega_v=d\tilde{v}$, we have that 
$g(\nabla_X^g v, X)=0$, and this tensorial map is such that
$I_v \neq \BOne_{\mc{H}}$.
Thus, if $\{ e_j\}_{i=0}^{2n}$ is a
local frame of eigenvectors of $C_v^t C_v$ with $e_0=v$, then     
$\nabla^g_{e_0}v=0$ globally, so $\lambda_0^2 =0$, and for $i>0$, the 
(local) vector field $\nabla^g_{e_i}v$ is horizontal and orthogonal to $e_i$, 
$\lambda_i^2=1$, $c_i^{-2}=2$, $(\nabla^g_{e_i}e_i)^{ver}=0$, and equation 
(\ref{eq10}) reduces to
$$
\frac{1}{2}\sum_{i> 0} \nabla^g_{e_i}\nabla^g_{e_i} v 
=-n v \, .  
$$

We let $X$ and $Y$ be locally defined horizontal vector fields  
(so orthogonal to $v$) nearby $p$. If $Y$
is orthogonal to both, $X$ and $I_vX=\nabla^g_Xv$, then so is 
$I_v Y= \nabla^g_Yv$.
Indeed, since $I_v$ preserves norms,
$$
g(I_vX,I_vY)=g(X,Y)= 0 \, ,
$$
and $I_vY$ is orthogonal to $I_vX$. On the other hand, 
we can complete $\{ X,I_vX\}$ to a normal frame 
$\{ X, I_vX, e_2, \ldots, e_{2n}\}$ of horizontal vectors such that 
$\nabla_{e_i}X\mid_p =0$. Then we have that $Y= \sum_j \beta_j e_j$, and 
therefore,
$$
g_p(X,I_v Y)=\sum \beta_i (p)g_p(X,I_v e_i)= -\sum \beta_i(p)
g_p(\nabla_{e_i}X,v)  =0 \, ,
$$
so $I_vY$ is orthogonal to $X$ also.  
It follows that the tensor $I_v$ leaves invariant any two plane 
horizontal distribution spanned by a set of the form $\{ X,\nabla_Xv\}$,
and as $I_v \neq \BOne$ is norm preserving, it must be that
$I_v^2 X= \nabla^g_{\nabla^g_Xv}v= -X$, and so
$I_v^2 = -\BOne_{\mc{H}}$. Notice that the projection 
of $\nabla^g_XY$ onto $v$ is tensorial in $X$ and $Y$, so we have that
$[X,Y]-[I_vX,I_vY] \in \mc{H}$, hence 
$[X,I_vY]+[I_vX,Y] \in \mc{H}$ also, 
and the triple $(\mb{S}^{2n+1},\mc{H},I_v)$
is an almost CR manifold. By defining $\varphi_v$ as $I_v$ on $\mc{H}$
and $\varphi_v(v)=0$, respectively, the quadruple $(v,\tilde{v},\varphi_v,g)$
is an almost contact metric structure, 
$$
g= d\tilde{v}\circ (\BOne \otimes \varphi_v) + \tilde{v}\otimes \tilde{v}\, ,  
$$        
and on horizontal vectors, 
$$
g(X,Y)= d\tilde{v}(X,I_vY)= \omega_v(X,I_vY) = -\tilde{v}([X,I_vY]) \, .
$$
 
We choose an orthonormal frame $\{e_i\}_{i=1}^{2n}$ of $\mc{H}$ adapted to
$I_v$, that is to say, $e_{2i}=I_v e_{2i-1}$, $i=1, \ldots, n$. The set  
$\{ v, e_1, \ldots, e_{2n}\}$ is a local frame of eigenvectors of
$C_v^t C_v$. If this frame is positively oriented, we extend the 
definition of $I_v$ to an almost complex structure $I$ to 
$T\mb{R}^{2n+2}\mid_{\mb{S}^{2n+1}}$ by declaring
that $I(r\partial_r\mid_{\mb{S}^{2n+1}})=v$, and then define a compatible
almost complex structure $I$ on the entire of $T\mb{R}^{2n+2}$ by the 
covariantly constant extension of the resulting tensor; if the said
frame is negatively oriented, we proceed similarly after declaring that 
$I( -r\partial_r\mid_{\mb{S}^{2n+1}})=v$. 
By construction, if $\{\partial_r,v, e_1, \ldots, e_{2n}\}$ is a positively 
oriented normal frame of $T\mb{R}^{2n+2}\mid_{\mb{S}^{2n+1}}$, we have that
$I\partial_r=v$, that the almost complex structure $I_v$ on $\mc{H}$ equals  
$I\mid_{\mc{H}}$, and that the structure $I$ is integrable if, and only if, 
$I_v$ is integrable. 

Consider the tensor $N \in C^{\infty}(\mb{S}^{2n+1},{\rm Hom}(\mc{H}\otimes
\mc{H}, \mc{H}))$ defined by   
$$
N(X,Y) = I_v([X,Y]-[I_vX,I_vY])-([I_vX,Y]+[X,I_vY]) \, .
$$
It is straightforward that $N(X,X)=0=N(X,I_v X)$, so $N \equiv 0$ if we prove
that $N(X,Y)=0$ for any $Y \in \mc{H}$ that is  orthogonal to both, $X$ and 
$I_v X$. Without loosing generality, we can assume that at the point $p$
where the computation is being carry out, $\nabla^g_u v \mid_p =0$ for
any pair $u,v \in \{ X, I_vX, Y, I_vY\}$. Then,
$$
\begin{array}{ccl}
{[X,Y]}_p & = & \nabla^g_X Y\mid_p - \nabla_Y^g X\mid_p = g(X,I_vY)v\mid_p-
g(Y,I_vX)v\mid_p =0 \, , \\  
{[I_vX,I_vY]}_p & = & \nabla^g_{I_vX} I_vY\mid_p - \nabla_{I_vY}^g I_vX\mid_p 
= g(I_vX,I^2_vY)v\mid_p- g(I_vY,I^2_vX)v\mid_p =0 \, ,  
\end{array}
$$
while
$$
\begin{array}{ccl}
{[I_vX,Y]}_p & = & \nabla^g_{I_vX} Y\mid_p - \nabla_Y^g I_vX\mid_p = 
g(I_vX,I_vY)v\mid_p-
g(Y,I^2_vX)v\mid_p = 0\, , \\  
{[X,I_vY]}_p & = & \nabla^g_{X} I_vY\mid_p - \nabla_{I_vY}^g X\mid_p 
= g(X,I^2_vY)v\mid_p- g(I_vY,I_vX)v\mid_p =0 \, .  
\end{array}
$$
Hence, $N(X,Y)=0$, as desired.
\qed

Since ${\rm dim}(\mb{S}\mb{O}(4))=6$, 
if $I,J,K$ denote the compatible complex structures in $\mb{R}^4=\mb{H}$ 
defined by the unit quaternions $i,j,k$, the subbundle of minimal unit fields
of $(\mb{S}^3,g)$ is given by
$$
\mc{V}_H(\mb{S}^3,g)=\{ v_H^{a,b,c}:=(aI+bJ+cK)(\partial_r\mid_{\mb{S}^3}) \, , 
\; a,b,c\in \mb{R}, \; a^2+b^2+c^1=1\} \, .  
$$
If $u\in C^{\infty}(\mb{S}^3; S^1(T\mb{S}^3))$, we define the {\it minimal
unit content} $v_{H}(u)$ of $u$ to be the $v_H^{a,b,c} \in \mc{V}_H(\mb{S}^3,g)$
for which the global inner product $\< u, v_H^{a,b,c}\>$ is the largest. Thus,
$$
v_{H}(u)= \frac{1}{C(u)} \left( \< u, I(\partial_r)\> I 
+ \< u, J(\partial_r)\> J 
+\< u, K(\partial_r)\> K\right)(\partial_r\mid_{\mb{S}^3}) \, , 
$$ 
where
$$
C(u)= \sqrt{ \< u, I(\partial_r)\>^2 + \< u, J(\partial_r)\>^2 + \< u, 
K(\partial_r)\>^2} \, .
$$ 
Notice that $f_{g(u)}(\mb{S}^3)$ and $f_{g(v_H(u))}(\mb{S}^3)$ are of the same
orientation, and lie in the same homology class in $S^1(T\mb{S}^3)$.

\begin{theorem}
Let $u$ be a unit vector field on $(\mb{S}^3, g)$, and let  
$v_{H}(u)$ be its minimal unit content. Then at any $(p,u(p)) \in
S^1(T\mb{S}^3)$, we have that
$$ 
d\mu_{f_{g(v_H(u))}}\mid_{T_{p,u(p)}f_{g(u)}(\mb{S}^3)} \leq 
d\mu_{f_{g(u)}}\mid_{p,u(p)} \, ,
$$ 
and the equality holds if, and only if, $u(p)=v_{H}(u)(p)$. Thus,  
$$
4\pi^2=\mu^{S^1(T\mb{S}^3)}_{v_H^{a,b,c}}(v_H^{a,b,c})=
\mu^{S^1(T\mb{S}^3)}_{g(v_H(u))}(v_H(u))\leq \mu^{S^1(T\mb{S}^3)}_{g(u)}(u)\, ,
$$
for any $v_H^{a,b,c}\in \mc{V}_H(\mb{S}^3,g)$, with equality on the right 
if, and only if, $u$ is minimal unit also.
\end{theorem} 

{\it Proof}. We consider any positively oriented normal frame 
(of eigenvectors of $C_u^tC_u$) of the form 
$\{ u_0:=u,u_1,u_2\}$, and let $\{ u_i^{f_{g(u)}}=c_i(u_i^{hor}+
(\nabla^g_{e_i} u)^{ver})\}_i$, $c_i^{2}=1+ \| \nabla^{g}_{u_i}u\|^2,$ be the
associated orthonormal tangent frame of $f_{g_{u}}(\mb{S}^3)$, and
$\{ \eta^i_{f_{g(u)}}\}$ be its associated coframe, respectively. Then, 
$d\mu_{f_{g(u)}(\mb{S}^3)}=\eta^0_{f_{g(u)}}\wedge \eta^1_{f_{g(u)}}\wedge  
\eta^2_{f_{g(u)}}$. If for convenience we denote by $I_{v_H(u)}$ the 
complex structure such that $v_H(u)=I_{v_H(u)}(\partial_r)$, we take  
also an orthonormal 
frame $\{ v_H(u),v_1,v_2:=I_{v_H(u)}(v_1)\}$ where 
$\nabla^g_{v_1}v_{H}(u)=I_{v_H(u)}(v_1)=v_2$
and $\nabla^g_{v_2}v_{H}(u)=I_{v_H(u)}(v_2)=-v_1$, with associated 
orthonormal tangent frame for $f_{g_{v_H(u)}}(\mb{S}^3)$ given by
$\{ v_0^{f_{g(v_H(u))}}=(v_{H}(u))^{hor},
v_1^{f_{g(v_H(u))}}=\frac{1}{\sqrt{2}}(v_1^{hor}+v_2^{ver}),
v_2^{f_{g(v_H(u))}}=\frac{1}{\sqrt{2}}(v_2^{hor}-v_1^{ver}) \}$, and let
$\{ \eta^i_{f_{g(v_H(u))}}\}$ be its associated coframe, so
$d\mu_{f_{g(v_H(u))}(\mb{S}^3)}=\eta^0_{f_{g(v_H(u))}}\wedge \eta^1_{f_{g(
v_H(u))}}\wedge  \eta^2_{f_{g(v_H(u))}}$. If 
$C^{v_H(u)}(u_{i}^{f_{g(u)}})$ is the vector of coordinates of 
$u_{i}^{f_{g(u)}}$ in the frame $\{v_{i}^{f_g(v_H(u))}\}$, and
$C^{v_H(u)}(\{ u_{i}^{f_{g(u)}}\})$ is the $3\times 3$ matrix
whose columns are $C^{v_H(u)}(u_{0}^{f_g(u)}), C^{v_H(u)}(u_{1}^{f_g(u)}),
C^{v_H(u)}(u_{2}^{f_{g(u)}})$, respectively, we have that  
$$
d\mu_{f_{g(v_H(u))}(\mb{S}^3)}(u_0^{f_{g(u)}},u_1^{f_{g(u)}},u_2^{f_{g(u)}})
=\det{C^{v_H(u)}(\{ u_{i}^{f_{g(u)}}\})} \leq 1 \, ,
$$
and the last equality occurs if, and only if, $u=v_H(u)$. 
This proves the first two assertions.  

We let $\iota: f_{g(u)}(\mb{S}^3) \hookrightarrow S^1(T\mb{S}^3)$ be the
inclusion map. Since $f_{g(u)}(\mb{S}^3)$ is homologous to 
$f_{g(v_H(u))}(\mb{S}^3)$ and the volume form is closed, the preliminary 
results imply that 
$$
\begin{array}{rcl}
\mu^{S^1(T\mb{S}^3)}_{g(v_H(u))}(v_H(u))={\displaystyle \int_{f_{g(v_H(u))}
(\mb{S}^3)} d\mu_{f_{g(v_H(u))}}} & = & {\displaystyle \int_{f_{g(u)}(\mb{S}^3)}
\iota^* d\mu_{f_{g(v_H(u))}}} \vspace{1mm}  \\ &  \leq & {\displaystyle   
\int_{f_{g(u)}(\mb{S}^3)} d\mu_{f_{g(u)}} = 
\mu^{S^1(T\mb{S}^3)}_{g(u)}(u)} \, .
\end{array}
$$
We finish the proof by observing that any $f_{g(v_H^{a,b,c})}(\mb{S}^3)$ 
is homologous to $f_{g(v_H(u))}(\mb{S}^3)$ and have 
$\mu^{S^1(T\mb{S}^3)}_{g(v_H^{a,b,c})}(v_H^{a,b,c})=\mu^{S^1(T\mb{S}^3)}_{
g(v_H(u))}(v_H(u))=2\mu_g(\mb{S}^3)$.
\qed
\medskip

We let $S_1^{2n+1}$ be a complete space form of curvature $1$. 
By the local nature of identity (\ref{eq10}), the
following result is straightforward, cf. \cite[Corollary 3.4]{ol}.  

\begin{example}
Under the projection map
$\pi: (\mb{S}^{2n+1},g) \rightarrow (S_1^{2n+1},\pi(g))$, any 
minimal unit vector field $v$ on $\mb{S}^{2n+1}$ passes to a minimal unit 
vector field on $S_1^{2n+1}$, and $\mu_{g(\pi (v))}^{TS_1^{2n+1}}(\pi(v))= 2^n 
\mu_{\pi(g)}(S_1^{2n+1})$.
\end{example}

\subsection{Product of spheres} 
For $r,s \in \mb{R}_{>0}$, $r^2+s^2=1$, and $k,l \in \mb{N}$,   
we consider the product manifold $(\mb{S}^k(r) \times \mb{S}^l(s),g_{k,l})$ 
isometrically embedded into $(\mb{S}^{k+l+1},
g)\hookrightarrow (\mb{R}^{k+l+2}, \| \, \|^2)$ as a hypersurface. 
We use the same complex structure always in the appropriate background 
Euclidean space leading to the definition of the Hopf vector fields on the 
factors, and denote by 
$v_r^k$ and $v_s^l$ their pullback under the corresponding projections to 
vector fields on the product. 

\begin{example}
If $v^k_r$ and $v^l_s$ are the standard Hopf vector fields of the factors of 
$(\mb{S}^k(r)\times \mb{S}^{l}(s),g_{k,l})$, respectively, 
for any $a,b\in \mb{R}$ such that $a^2+b^2=1$, we let
$v_{a,b}:=av^k_r+bv^l_s$. Then:
\begin{enumerate}[label={\rm (\alph*)}]
\item If $ab=0$, where defined, $v_{a,b}/\| v_{a,b}\|$ is minimal unit,
and when at least one of $k$ or $l$ is even, no other normalized $v_{a,b}$ is 
minimal unit. 
\item If both $k$ and $l$ are odd, then $\| v_{\frac{1}{\sqrt{2}},\pm 
\frac{1}{ \sqrt{2}}}\|=1$, and  $\pm \frac{1}{\sqrt{2}}v_{1,\pm 1}$ is 
minimal unit if, and only if, $k=l=1$, or $k=l$ and $r=s$. 
\end{enumerate}
\end{example}

This follows readily by the arguments of Example \ref{hfs}, and can 
be expanded to say that when $ab=0$, and both $k$ and $l$ are odd integers, 
or when exactly one of $k$ and $l$ is even and its associated coefficient 
vanishes, $f_{g}(tv_{a,b}) (M)$ is a canonical cycle in its homology class 
which converges uniformly to the canonical zero section cycle as 
$t\rightarrow 0$. This remains so for the minimal unit cases in (b) also. 
Notice that when $k$ and $l$ are both odd, 
$\pm v_{1,\pm 1}$ is a Killing field no matter the values of
$r$ and $s$, and so (b) shows that a unit Killing 
field does not have to be minimal.     

\begin{remark}
The Calabi-Eckmann manifolds above where one of the factors is a circle 
exemplify the class of compact locally conformally K\"ahler manifolds that 
carry a nontrivial parallel field. These manifolds must be either locally 
conformally K\"ahler with a parallel Lee form, or globally conformally 
K\"ahler equal to a suitable quotient of the flat plane $\mb{R}^{2}$ times a 
conformal deformation of a K\"ahler manifold of dimension $k+l-2>2$ 
\cite{moro}. 
The normalized parallel vector field on these manifolds is minimal unit, and
its foliation is tautologically taut.    
\end{remark}

\subsection{The angular vector field on uniformized Riemann surfaces}
A closed oriented Riemannian surface $\Sigma=\Sigma^k$ of genus $k=0$, $k=1$ 
or $k>1$ 
and metric of constant scalar curvature $2$, $0$ or $-2$, is 
$\mb{S}^2$, $\mb{C}/(\mb{Z}+ \tau \mb{Z})$ for some specific 
$\tau=\tau_{\Sigma}$,  or $D/\Gamma$, $D$ the Poincar\'e disk $D$ and 
$\Gamma=\Gamma_{\Sigma}$ a specific 
Fuchsian group, respectively. If $\mc{U}_{\Sigma}$ denotes the universal
cover of $\Sigma$ with its canonical metric, and $\Gamma_{\Sigma}$ the 
discrete group that acts on it by isometries to defines $\Sigma$ as a 
quotient, we let 
$$
\pi_{\Sigma}: \mc{U}_{\Sigma} \rightarrow 
\Sigma=\mc{U}_{\Sigma}/\Gamma_{\Sigma}
$$
denote the Riemannian submersion that is obtained. 

In geodesic polar coordinates, the lift of the metric on $\Sigma$ to the 
$\mc{U}_{\Sigma}$ is of the form $dr^2+ \sin^2{r}\, d\theta^2$, 
$dr^2 + r^2 \, d\theta^2$ or $dr^2 + \sinh^2{r} \, d\theta^2$, respectively.
We define the angular vector field $e_{\theta}$ on the cover to be
$\sin{r}\partial_{\theta}$, $r\partial_{\theta}$, or  
$\sinh{r}\partial_{\theta}$, according to the case. It is a smooth vector
field with simple zeroes at antipodal points determined by the equation  
$\sin{r}=0$ when $k=0$, or a simple zero at $r=0$ in any of the other cases.
 
\begin{example}
The field $e_{\theta}$ on $\mc{U}_{\Sigma}$ is a minimal singular unit vector 
field, and under the projection map $\pi_{\Sigma}$,  
it passes to a minimal singular unit vector field on $\Sigma$. 
\end{example}

The verification of the identity (\ref{eq10}) for $e_{\theta}/\| e_{\theta}\|$
is straightforward (cf. Example \ref{hfs}(c) when $2n=2$). This local 
condition then holds for the 
field $\pi_{\Sigma \; *}(e_{\theta})$, whose additional singularities are     
determined by the identifications that the group $\Gamma_{\Sigma}$
introduces. In any of the cases, the field $e_r=r\partial_r$ is minimal
singular unit orthogonal to $e_{\theta}$, and 
$\pi_{\Sigma \; *}(e_r)$ is minimal singular unit on $\Sigma$,
orthogonal to $\pi_{\Sigma \; *}(e_{\theta})$ also.

\end{document}